\documentclass[a4paper,12pt,reqno]{amsart}
\usepackage{amssymb}
\usepackage[all,cmtip]{xy}
\usepackage{amsmath}
\usepackage{chemarrow}
\usepackage{extarrows}
\usepackage{xcolor}
\usepackage[colorlinks=true, linkcolor=blue!55!black, citecolor=green!40!black, urlcolor=blue!65!black]{hyperref}
\usepackage{cleveref}
\usepackage{mathrsfs}
\usepackage{titletoc}
\usepackage{bm}
\usepackage{dsfont} 
\usepackage{geometry}
\usepackage{colordvi}
\usepackage{color}
\usepackage{stmaryrd}
\usepackage{enumitem}
\usepackage{tikz-cd}

\allowdisplaybreaks

\usepackage{calrsfs}
\DeclareMathAlphabet{\cal}{OMS}{zplm}{m}{n}

\DeclareMathOperator{\ad}{\mathrm{ad}}

\DeclareMathOperator{\Aut}{Aut}

\DeclareMathOperator{\cC}{{\cal C}}

\DeclareMathOperator{\cF}{{\cal F}}
\DeclareMathOperator{\cG}{{\cal G}}
\DeclareMathOperator{\cha}{\mathrm{char}}
\DeclareMathOperator{\cM}{{\cal M}}

\DeclareMathOperator{\Coker}{Coker}
\DeclareMathOperator{\cop}{cop}
\DeclareMathOperator{\cS}{{\cal S}}
\DeclareMathOperator{\cZ}{{\cal Z}}

\DeclareMathOperator{\Ext}{Ext}

\DeclareMathOperator{\Hom}{Hom}

\DeclareMathOperator{\id}{id}

\DeclareMathOperator{\Ker}{Ker}

\DeclareMathOperator{\lhu}{\leftharpoonup}

\DeclareMathOperator{\Max}{Maxspec}

\DeclareMathOperator{\op}{op}

\DeclareMathOperator{\rhu}{\rightharpoonup}

\DeclareMathOperator{\vep}{\varepsilon}

\DeclareMathOperator{\YD}{\cal{YD}}
\DeclareMathOperator{\Z}{\mathbb{Z}}

\newcommand{\YDlr}{{}_{H}{\cal YD}^{H}}   
\newcommand{\YDrr}{{\cal YD}^{H}_{H}}     
\newcommand{\YDrl}{{}^{H}{\cal YD}_{H}}   
\newcommand{\YDll}{{}^{H}_{H}{\cal YD}}   

\newcommand{\lint}{\textstyle \int_{H}^{\ell}}
\newcommand{\rint}{\textstyle \int_{H}^{r}}
\newcommand{\lcoint}{\textstyle \int^{H}_{\ell}}
\newcommand{\rcoint}{\textstyle \int^{H}_{r}}

\newcommand{\To}{\longrightarrow}

\numberwithin{equation}{section}

\theoremstyle{definition}

\newtheorem{thm}{Theorem}[section]
\newtheorem{prop}[thm]{Proposition}
\newtheorem{lem}[thm]{Lemma}

\newtheorem{cor}[thm]{Corollary}

\newtheorem{defn}[thm]{Definition}
\newtheorem{rem}[thm]{Remark}
\newtheorem{ques}[thm]{Question}
\newtheorem{ex}[thm]{Example}

\Crefname{defn}{Definition}{Definitions}
\Crefname{lem}{Lemma}{Lemmas}
\Crefname{prop}{Proposition}{Propositions}
\Crefname{cor}{Corollary}{Corollaries}
\Crefname{thm}{Theorem}{Theorems}
\Crefname{rem}{Remark}{Remarks}
\Crefname{ques}{Question}{Questions}
\Crefname{ex}{Example}{Examples}

\hypersetup{pdftitle={The Fourth-Power Formula for the Antipode in AS-Gorenstein Hopf Algebras}, pdfauthor={Quanshui Wu and Ruipeng Zhu}}

\begin{document}
	
\title[The Fourth-Power Formula for the Antipode in AS-Gorenstein Hopf Algebras]{The Fourth-Power Formula for the Antipode\\ in AS-Gorenstein Hopf Algebras}
\author{Quanshui Wu}
\address{School of Mathematical Sciences, Fudan University, Shanghai 200433, China}
\email{qswu@fudan.edu.cn}

\author{Ruipeng Zhu}
\address{School of Mathematics, Shanghai University of Finance and Economics, Shanghai 200433, China}
\email{zhuruipeng@sufe.edu.cn}

\begin{abstract}
	AS-Gorenstein Hopf algebras, for which homological integrals are defined, form a broad class of (infinite-dimensional) Hopf algebras. We establish an $S^4$ formula for left AS-Gorenstein Hopf algebras, thereby generalizing Radford's $S^4$ formula for finite-dimensional Hopf algebras and refining the formula of Brown and Zhang by identifying a group-like conjugating element in the Hopf algebra $H$.
  Our approach is to endow the space of left homological integrals $\lint$ with a canonical right $H^*$-module structure compatible with its natural right $H$-action.
  If the $H^*$-action is rational, then $\lint$ carries a canonical $(\id_H,S^4)$-twisted right-left Yetter--Drinfeld module structure.
  The right $H$-action and left $H$-coaction on this one-dimensional space determine a character and a group-like element in $H$, respectively; the twisted Yetter--Drinfeld compatibility condition then yields the $S^4$ formula. The $H^*$-action on $\lint$ is rational if either the coradical of $H$ is coreflexive or the trivial left module $\Bbbk$ admits a finitely generated projective resolution.
  In particular, the latter condition holds when $H$ is left noetherian.
    We also establish a corresponding result for right AS-co-Gorenstein Hopf algebras with coreflexive coradical. 
    As applications, we use the $S^4$ formula to study the bijectivity of the antipode and derive explicit finite-order bounds when $H$ is module-finite and faithfully flat over a central right coideal subalgebra, or finitely cogenerated and faithfully coflat over a cocentral quotient.

\end{abstract}


\thanks{This research has been supported by the NSFC (Grant Nos. 12471032 and 12301052).}

\maketitle

\section{Introduction}

A key 
result in the theory of finite-dimensional Hopf algebras $H$ is Radford's $S^4$ formula, which expresses the fourth power of the antipode $S$ in terms of the distinguished group-like elements $g \in G(H)$ and $\alpha \in G(H^*)$ arising from the right integral in $H^*$ and the left integral in $H$, respectively: for all $h\in H$,
\begin{equation}\label{eq:Radford}
	S^4(h) = g \, (\alpha \rhu h \lhu \alpha^{-1}) \, g^{-1}.
\end{equation}
The formula was first proved by Larson \cite[Theorem~5.5]{Lar71} for finite-dimensional unimodular Hopf algebras and was extended by Radford \cite[Proposition~6]{Rad76} to arbitrary finite-dimensional Hopf algebras.
Schneider \cite[Theorem~3.8]{Sch95} later provided a more conceptual proof based on the Frobenius structure of $H$, deriving the formula by comparing the two expressions for the Nakayama automorphism obtained from the integrals of $H$ and $H^{*}$, respectively.

Since then, Radford's $S^{4}$-formula has been generalized to
finite-dimensional quasi-Hopf algebras \cite[Corollary~6.3]{HN99},
finite-dimensional braided Hopf algebras \cite[Theorem~3.6]{BKLT00}, 
finite-dimensional weak Hopf algebras \cite[Theorem~5.13]{Nik02}, 
co-Frobenius Hopf algebras \cite[Theorem~2.8]{BBT07}, and
finite tensor categories \cite[Theorem~3.3]{ENO04}.


The study of the homological properties of (noetherian) infinite-dimensional Hopf algebras has become a central theme in noncommutative algebra and quantum group theory.
In their seminal work \cite{AS87}, Artin and Schelter introduced what are now called Artin--Schelter regular algebras (AS-regular, for short) as noncommutative analogues of polynomial rings, thereby initiating the modern development of noncommutative projective geometry.
The Artin--Schelter Gorenstein (AS-Gorenstein, for short) condition subsequently emerged as a natural generalization, retaining the homological duality structure while relaxing the hypothesis of finite global dimension.

The ring-theoretic and homological study of noetherian Hopf algebras dates back to the work of Brown and Goodearl \cite{BG97, Bro98}; see \cite{Bro26} for a survey.
A central question in this area is:
 \centerline{\it Is every noetherian Hopf algebra AS-Gorenstein?}\\
This question is known as the Brown--Goodearl conjecture. 
The conjecture was originally formulated for noetherian affine Hopf algebras satisfying a polynomial identity and was proved in this setting by the first author and Zhang \cite{WZ03}. Consequently, many quantum groups at roots of unity are AS-Gorenstein.
The class of AS-Gorenstein Hopf algebras is broad 
\cite{BG97, Bro98, BZ08, Bro26}.
Other examples include universal enveloping algebras of finite-dimensional Lie algebras and group algebras of polycyclic-by-finite groups \cite[Sections~6.3 and~6.7]{BZ08}.
Quantum groups also provide important non-PI examples.
Let $\mathfrak g$ be a finite-dimensional complex semisimple Lie algebra, let $G$ be a connected complex semisimple algebraic group, and let $q$ be an indeterminate.
Then the Drinfeld--Jimbo quantized enveloping algebra $U_q(\mathfrak g)$ and the standard quantized coordinate algebra ${\cal O}_q(G)$ over $\mathbb C(q)$ are both noetherian AS-regular, and hence AS-Gorenstein, with global dimensions $\dim_{\mathbb C}\mathfrak g$ and $\dim G$, respectively; see \cite[Sections~6.4--6.5]{BZ08}.

The AS-Gorenstein property (\Cref{def:ASG}) is closely connected to the symmetry structure of Hopf algebras through the notion of homological integrals introduced by Lu, the first author, and Zhang \cite{LWZ07}.
This invariant links the homological properties of Hopf algebras to their representation theory and encodes information about the rigid dualizing complex, the Nakayama automorphism \cite{BZ08}, and the Calabi--Yau property \cite{HVOZ10}.
It has also played an important role in the classification of Hopf algebras of low Gelfand--Kirillov dimension \cite{BZ10, WLD16, Liu20}.

Brown and Zhang \cite{BZ08} extended Radford's formula for the fourth power of the antipode from finite-dimensional Hopf algebras to noetherian AS-Gorenstein Hopf algebras.
In the spirit of Schneider's approach, their proof compares two computations of the Nakayama automorphisms of the AS-Gorenstein Hopf algebra $H$. 
Let $H$ be such an algebra, and let $\alpha \colon H \to \Bbbk$ denote the character determined by the right $H$-action on the space of left homological integrals.
The generalized Radford formula then takes the form
\begin{equation}\label{eq:Rad-Bro-Zha}
	S^4(h)
	\, = \,
	u \, (\alpha \rhu h \lhu \alpha^{-1}) \, u^{-1},
\end{equation}
where $u$ is an invertible element of $H$.
When $H$ is finite-dimensional, $\alpha$ coincides with the distinguished group-like element in $G(H^{*})$, and $u$ can be taken to be the distinguished group-like element of $H$. 
This result provides a homological interpretation of the fourth power of the antipode and has become a fundamental tool in the study of Nakayama automorphisms, homological determinants, and twisted Calabi--Yau properties of Hopf algebras.

The following question was first posed by Brown and Zhang \cite[Question~4.6]{BZ08}; see also \cite[Question~3.2]{LOWY18} and \cite[Question~4.4]{Bro26}.

\begin{ques}
	Which inner automorphism of $H$ is induced by conjugation by the element $u$ appearing in \eqref{eq:Rad-Bro-Zha}?
\end{ques}

We address this question by constructing a canonical group-like element from the space of homological integrals, thereby obtaining an $S^4$ formula of the same shape as \eqref{eq:Radford} for AS-Gorenstein Hopf algebras.
A key structural property is the existence of a natural right $H^*$-action on $\Ext_H^i({}_{\vep}\Bbbk,H)$ that is compatible with the usual right $H$-action through an $S^4$-twisted identity.
These actions are constructed using suitable resolutions in Yetter--Drinfeld categories and are independent of the choice of resolution. The construction also applies to more general Ext groups (\Cref{prop:Ext-structure}).

The essential issue in passing from these actions to a twisted Yetter--Drinfeld module structure (\Cref{defn:twisted-YD}) is rationality.
A right $H^*$-module is rational if its action is induced by a left $H$-coaction.
For infinite-dimensional $H$, the one-dimensionality of an $H^*$-module alone does not imply rationality.
We prove that the canonical $H^*$-action on $\lint$ is rational under either of two hypotheses: coreflexivity of the coradical of $H$ or a suitable finiteness condition on the trivial $H$-module.
Both hypotheses ensure that the canonical actions on $\lint$ determine an $(\id_H,S^4)$-twisted right-left Yetter--Drinfeld module structure.
The resulting coaction on this one-dimensional space determines the desired group-like element.

To state the theorem precisely, we recall the relevant definitions.
A coalgebra $C$ is \emph{coreflexive} if the canonical embedding $C \hookrightarrow (C^{*})^{\circ}$ is an isomorphism.
A module $M$ over a ring $R$ is of type $\mathrm{FP}_{d}$ for some $d \geqslant 0$ if it admits a projective resolution
\[ \cdots \To P_{d+1} \To P_d \To \cdots \To P_1 \To P_0 \To M \To 0 \]
in which $P_0, P_1, \dots, P_d$ are finitely generated.
If $M$ admits a projective resolution in which each $P_i$ is finitely generated, then $M$ is said to be of type $\mathrm{FP}_\infty$ (also called pseudo-coherent in some literature).

\begin{thm}[\Cref{thm:main}]\label{thm:intro-main}
	Let $H$ be a left AS-Gorenstein Hopf algebra of dimension $d$ with antipode $S$.
	Suppose that either the coradical of $H$ is coreflexive or the trivial left $H$-module ${}_{\vep}\Bbbk$ is of type $\mathrm{FP}_{d}$.
	Then there exist a group-like element $g \in G(H)$ and a character $\alpha \colon H \to \Bbbk$, both determined by the space of left homological integrals of $H$, such that, for all $h \in H$,
  \begin{equation}\label{eq:S4}
    S^4(h) = g \, \bigl( \alpha \rhu h \lhu \alpha^{-1} \bigr) \, g^{-1}.
  \end{equation}
\end{thm}

For a noetherian AS-Gorenstein Hopf algebra $H$ of dimension $d$, the trivial left $H$-module ${}_{\vep}\Bbbk$ is of type $\mathrm{FP}_{d}$.
Consequently, the $S^{4}$ formula in \Cref{thm:intro-main} applies in this case, refining the inner-automorphism version of \eqref{eq:Rad-Bro-Zha}.
In contrast to the approach of Brown and Zhang \cite{BZ08}, which relies on the Nakayama automorphism, we directly exhibit a canonical twisted Yetter--Drinfeld module structure (\Cref{defn:twisted-YD}) on the space of homological integrals (\Cref{thm:main}). The $S^{4}$ formula then follows from an elementary computation using the twisted Yetter--Drinfeld compatibility condition.
The same argument yields the corresponding $S^4$ formula for a right AS-co-Gorenstein Hopf algebra whose right homological cointegral space is a rational left $H^*$-module, generalizing the classical formula for co-Frobenius Hopf algebras. 
If $H$ is both left AS-Gorenstein and right AS-co-Gorenstein and has a coreflexive coradical, then the spaces of homological integrals and homological cointegrals determine pairs $(g, \alpha)$ and $(\widetilde{g}, \widetilde{\alpha})$, respectively, each consisting of a group-like element and a character.
These pairs need not coincide; see \Cref{ex:klein}.

\Cref{sec:app} presents applications of the $S^4$ formulas established in \Cref{thm:main} and \Cref{thm:main-co}.
In particular, every left AS-Gorenstein Hopf algebra satisfying the hypothesis of \Cref{thm:intro-main} has a bijective antipode (\Cref{prop:bijective-antipode-ASG}(1)), and every right AS-co-Gorenstein Hopf algebra with coreflexive coradical has a bijective antipode (\Cref{prop:bijective-antipode-ASG}(2)).
We then apply the $S^4$ formulas to address the question posed by Brown and Zhang \cite[Question~6.2]{BZ08}: Does the antipode of a noetherian affine PI Hopf algebra have finite order? 
To this end, let $H$ be a Hopf algebra that is a finite module over a central Hopf subalgebra. If $H$ is AS-Gorenstein of dimension $d$ and either its coradical is coreflexive or ${}_{\vep}\Bbbk$ is of type $\mathrm{FP}_{d}$, then $S^2$ has finite order (\Cref{cor:order-ASG}). If $H$ is AS-co-Gorenstein and its coradical is coreflexive, then $S^2$ has finite order (\Cref{prop:order-AScoG-B}).
These results provide partial answers to Brown and Zhang's question.

\medskip
\noindent\textbf{Organization.}
\Cref{sec:prelim} reviews the necessary background on AS-Gorenstein and AS-co-Gorenstein Hopf algebras, homological integrals and cointegrals, twisted Yetter--Drinfeld modules, and twisted Drinfeld centers.
In \Cref{sec:S4}, we construct a twisted Yetter--Drinfeld structure on the space of homological integrals and prove the $S^{4}$ formula for AS-Gorenstein and AS-co-Gorenstein Hopf algebras.
\Cref{sec:app} presents examples and applications, including the bijectivity of the antipode and the order of $S^{2}$.

\medskip
\noindent\textbf{Conventions and notation.}
Throughout, $\Bbbk$ is a field, and all unadorned tensor products and $\Hom$ are taken over $\Bbbk$.
For a coalgebra $C$, we use Sweedler's notation $\Delta(c) = c_1 \otimes c_2$ (with summation understood). Similarly, we write a left coaction as $\rho(v) = v_{-1} \otimes v_0$ and a right coaction as $\rho(v) = v_0\otimes v_1$.
For a Hopf algebra $H$, we denote its unit map by $\eta$, its counit map by $\vep$, its antipode by $S$, and its group of group-like elements by $G(H)$.
For $h, l \in H$ and $h^* \in H^* := \Hom(H, \Bbbk)$, we write $h \lhu h^* = \langle h^*, h_1 \rangle h_2$, $h^* \rhu h = h_1 \langle h^*, h_2 \rangle$, $(h \rhu h^*)(l) = \langle h^*, lh \rangle$, and $(h^* \lhu h)(l) = \langle h^*, hl \rangle$.
Let $\varphi \colon A \to B$ be an algebra homomorphism, and let $M$ be a left $B$-module.
We denote by ${}_{\varphi}M$ the left $A$-module obtained from $M$ by restriction of scalars along $\varphi$.
Dually, let $\psi \colon C \to D$ be a coalgebra homomorphism, and let $N$ be a left $C$-comodule with coaction $\rho \colon N \to C \otimes N$.
We denote by ${}^{\psi}N$ the left $D$-comodule obtained from $N$ by corestriction of scalars along $\psi$, whose coaction is given by the composite $N \xrightarrow{\rho} C \otimes N \xrightarrow{\psi \otimes \id_N} D \otimes N$.
The right module and right comodule cases are denoted similarly.

\section{Preliminaries}\label{sec:prelim}

The $S^4$ formula is stated in terms of homological integrals and twisted Yetter--Drinfeld modules.
We fix these definitions here, starting with AS-Gorenstein Hopf algebras.

\subsection{AS-Gorenstein Hopf algebras and their homological integrals}\label{subsec:ASG}

We first recall the Artin--Schelter conditions for Hopf algebras; see \cite{WZ03, BZ08}.

\begin{defn}\label{def:ASG}
	A Hopf algebra $H$ is called \emph{left Artin--Schelter Gorenstein} (left AS-Gorenstein, for short) of dimension $d$ if
	\begin{enumerate}[label=\textup{(\arabic*)},leftmargin=2.2em]
		\item $_HH$ has finite injective dimension $d$;
		\item $\Ext^i_H({}_{\vep}\Bbbk, H) \cong \begin{cases}
			0, & i \neq d,\\
			\Bbbk, & i = d,
		\end{cases}$
		as vector spaces, where ${}_{\vep}\Bbbk$ denotes the trivial left $H$-module, on which $H$ acts via the counit $\vep$.
	\end{enumerate}
 \end{defn}   
	A right AS-Gorenstein Hopf algebra is defined analogously. A Hopf algebra is called \emph{AS-Gorenstein} if it is both left and right AS-Gorenstein.
	If an AS-Gorenstein Hopf algebra $H$ has finite global dimension $d$, then $H$ is called \emph{AS-regular} of dimension $d$.
	The definition of AS-regularity used here does not require the algebra to have finite Gelfand--Kirillov dimension.

The Lorenz--Lorenz lemma \cite{LL95} states that the left global dimension of a Hopf algebra $H$ equals the projective dimension of 
${}_{\vep}\Bbbk$.

Let $H$ be a left AS-Gorenstein Hopf algebra of dimension $d$.
The one-dimensional space
\[ \lint := \Ext^d_H({}_{\vep}\Bbbk, H) \]
is called the \emph{space of left homological integrals} of $H$, and any nonzero element of $\lint$ is called a \emph{left homological integral} \cite{LWZ07}.
The space $\lint$ has a natural right $H$-module structure induced by right multiplication in the second argument, and there is a unique algebra homomorphism $\alpha \colon H \to \Bbbk$ such that $\lint \cong \Bbbk_{\alpha}$ as right $H$-modules.
For a right AS-Gorenstein Hopf algebra $H$, the \emph{space of right homological integrals} is defined analogously by $\rint := \Ext^d_{H^{\op}}(\Bbbk_{\vep}, H)$.
When $H$ is finite-dimensional, it is AS-Gorenstein of dimension $0$, and its homological integral coincides with the usual integral \cite{LWZ07}.


In contrast to the definition in \cite{BZ08}, \Cref{def:ASG} does not require $H$ to be noetherian.
Many known AS-regular Hopf algebras satisfy a weaker finiteness condition: the trivial module ${}_{\vep}\Bbbk$ is of type $\mathrm{FP}_{\infty}$.
A group $G$ is called a Poincar\'e duality group over $\Bbbk$ if the group algebra $\Bbbk G$ is AS-regular and the trivial $\Bbbk G$-module $\Bbbk$ is of type $\mathrm{FP}_{\infty}$.
Group algebras of Poincar\'e duality groups provide a rich source of non-noetherian AS-regular Hopf algebras.
Let $X$ be a closed, connected, oriented, aspherical smooth manifold of dimension $d$, and let $G = \pi_1(X)$ be the fundamental group of $X$. 
Then $G$ is a Poincar\'e duality group over $\mathbb{C}$ of dimension $d$ with trivial orientation character \cite[Ch.~VIII]{Bro82}.
Consequently, the group algebra $\mathbb{C}G$ is AS-regular of dimension $d$ and is, in fact, a Calabi--Yau algebra of dimension $d$ in the sense of Ginzburg \cite{Gin07}.
Let $\Sigma_g$ be a closed orientable surface of genus $g \geqslant 2$ and set $G = \pi_1(\Sigma_g)$.
Then $\mathbb{C}G$ is AS-regular of dimension $2$.
Since the commutator subgroup $[G, G]$ is a free group of infinite rank \cite[Theorem~B]{Put22}, the group algebra $\mathbb{C}G$ is not noetherian.

\subsection{AS-co-Gorenstein Hopf algebras and their homological cointegrals}\label{subsec:AScoG}

The category of right comodules over a coalgebra is a Grothendieck category, but it may not have enough projective objects.
Therefore, the right global dimension of a coalgebra $C$ is defined as the supremum of the injective dimensions of all right $C$-comodules \cite{NTZ96}.

\begin{defn}\label{def:AScoG}
	A Hopf algebra $H$ is called \emph{right AS-co-Gorenstein} of dimension $d$ if
	\begin{enumerate}[label=\textup{(\arabic*)},leftmargin=2.2em]
		\item $\Ext^i_{\cM^H}(H, N) = 0$ for any right $H$-comodule $N$ and $i > d$;
		\item $\Ext^i_{\cM^H}(H, \Bbbk^{\eta}) \cong
		\begin{cases}
			0, & i \neq d,\\
			\Bbbk, & i = d,
		\end{cases}$
		as vector spaces, where $\Bbbk^{\eta}$ denotes the trivial right $H$-comodule, whose coaction is induced by the unit $\eta$.
	\end{enumerate}
    \end{defn}
    
A left AS-co-Gorenstein Hopf algebra is defined analogously. A Hopf algebra is called \emph{AS-co-Gorenstein} if it is both left and right AS-co-Gorenstein.
If an AS-co-Gorenstein Hopf algebra $H$ has finite global dimension $d$ as a coalgebra, then $H$ is called \emph{AS-co-regular} of dimension $d$.

Analogous to the Lorenz--Lorenz lemma \cite{LL95}, the right global dimension of a Hopf algebra $H$ as a coalgebra equals the injective dimension of the trivial right $H$-comodule $\Bbbk^{\eta}$; see, for example, \cite[Theorem~1.1]{DNT01a}.

\begin{rem}
	A coalgebra $C$ is called \emph{AS-regular} of dimension $d$ \cite[Remark~2.2]{HTOZ14} if the following conditions hold: $C$ has finite right global dimension; for every simple right $C$-comodule $T$, $\Ext_{\cM^C}^i(C, T) = 0$ for all $i \neq d$ and $\Ext_{\cM^C}^d(C, T)$ is a simple right $C$-comodule; and the analogous conditions hold for left $C$-comodules.
	
	Let $H$ be a Hopf algebra with a bijective antipode.
    A routine verification shows that $H$ is AS-co-regular of dimension $d$ if and only if $H$ is an AS-regular coalgebra of dimension $d$.
	If $H^*$ is noetherian as an algebra, He, Torrecillas, Van Oystaeyen and Zhang \cite[Proposition~2.3]{HTOZ14} proved that $H$ is an AS-regular coalgebra if and only if $H^*$ is an AS-regular algebra.
	Let ${\cal U}(\mathfrak{g})$ be the enveloping algebra of a finite-dimensional Lie algebra  $\mathfrak{g}$, and assume that $\cha(\Bbbk) = 0$.
	Then ${\cal U}(\mathfrak{g})^* \cong \Bbbk[\![x_1, x_2, \dots, x_n]\!]$ as algebras, where $n = \dim_{\Bbbk} \mathfrak{g}$. The latter is a noetherian AS-regular algebra of dimension $n$.
	Consequently, ${\cal U}(\mathfrak{g})$ is AS-co-regular of dimension $n$.
\end{rem}

Let $H$ be a right AS-co-Gorenstein Hopf algebra of dimension $d$.
The one-dimensional space
\[ \rcoint := \Ext^d_{\cM^H}(H, \Bbbk^{\eta}) \]
is called the \emph{right homological cointegral space} of $H$, and any nonzero element of $\rcoint$ is called a \emph{right homological cointegral}.
Moreover, $\rcoint$ has a natural left $H^*$-module structure induced by the left $H$-comodule structure on $H$ in the first variable (see \Cref{lem:Hom-structure-co} and \Cref{thm:main-co}).
For a left AS-co-Gorenstein Hopf algebra $H$, the \emph{left homological cointegral space} is defined analogously by $\lcoint := \Ext^d_{{}^H\!\!\cM}(H, {}^{\eta}\Bbbk)$.

A Hopf algebra $H$ is AS-co-Gorenstein of dimension $0$ if and only if it is co-Frobenius; that is, there exists a nonzero element $\lambda \in H^*$ such that $\lambda * f = \langle f, 1_H \rangle \lambda$ for any $f \in H^*$.
In this case, any such nonzero element is a right homological cointegral, as defined above.

\subsection{Twisted Yetter--Drinfeld modules}\label{subsec:twisted-YD}

The notion of Yetter--Drinfeld modules was introduced by Yetter \cite{Yet90}, who termed them crossed bimodules; the present name was subsequently coined by Radford and Towber in \cite{RT93}.
A Yetter--Drinfeld module structure naturally gives a solution to the Yang--Baxter equation \cite[Proposition 3]{RT93}.
Their categorical significance stems from the equivalence between the category of Yetter--Drinfeld modules over a Hopf algebra $H$ and the Drinfeld center of the category of $H$-(co)modules \cite{Kas95, EGNO15}. 
A twisted notion of Yetter--Drinfeld modules was introduced by Panaite and Staic \cite{PS07}. 
In this subsection, we give a variant of the twisted notion of Panaite and Staic \cite{PS07} in which the Yetter--Drinfeld compatibility condition is twisted by two Hopf algebra endomorphisms rather than by automorphisms.

Let $H$ be a $\Bbbk$-bialgebra, and let $\sigma$ and $\tau$ be two bialgebra endomorphisms of $H$.

\begin{defn}\cite[Definition~2.1]{PS07}
	\label{defn:twisted-YD}
	A \emph{$(\sigma,\tau)$-twisted left-right Yetter--Drinfeld module} is a triple $(V, \rhd, \rho)$ consisting of a vector space $V$ endowed with a left $H$-module structure $(V, \rhd)$ and a right $H$-comodule structure $(V, \rho)$, 
	such that for all $v \in V$ and $h \in H$,
	\begin{equation}\label{eq:yd-lr}
		h_1 \rhd v_0 \otimes \tau(h_2) v_1
		\;=\;
		(h_2 \rhd v)_0 \otimes (h_2 \rhd v)_1 \sigma(h_1).
	\end{equation}
	We denote by $\YDlr(\sigma, \tau)$ the category of $(\sigma, \tau)$-twisted left-right Yetter--Drinfeld modules, with morphisms given by the maps that are both $H$-linear and $H$-colinear.
\end{defn}

Similarly, the categories $\YDrr(\sigma, \tau)$, $\YDrl(\sigma, \tau)$, and
$\YDll(\sigma, \tau)$ of \emph{$(\sigma, \tau)$-twisted right-right, right-left, and left-left Yetter--Drinfeld modules}, respectively, are defined by replacing \eqref{eq:yd-lr} with the compatibility conditions
\begin{align}
	v_0 \lhd h_1 \otimes v_1\sigma(h_2) & = (v \lhd h_2)_0 \otimes \tau(h_1) (v \lhd h_2)_1, \label{eq:yd-rr} \\
	v_{-1}\sigma(h_1) \otimes v_0 \lhd h_2 & = \tau(h_2) (v \lhd h_1)_{-1} \otimes (v \lhd h_1)_0, \label{eq:yd-rl} \\
	\tau(h_1)v_{-1} \otimes h_2 \rhd v_0 & = (h_1 \rhd v)_{-1} \sigma(h_2) \otimes (h_1 \rhd v)_0. \label{eq:yd-ll}
\end{align}
When $\sigma = \tau = \id_H$, these reduce to the classical Yetter--Drinfeld categories of \cite{RT93}; for example, $\YDlr(\id_H,\id_H)$ is the usual category $\YDlr$.

A right (resp.\ left) $H$-comodule structure on $V$ induces a left (resp.\ right) $H^*$-module structure via $h^* \rhd v := v_0 \langle h^*, v_1 \rangle$ (resp.\ $v \lhd h^* := \langle h^*, v_{-1} \rangle v_0$).
Applying $\langle h^*, - \rangle$, $h^* \in H^*$, to the comodule tensorand of \eqref{eq:yd-lr}--\eqref{eq:yd-ll} shows that these conditions are equivalent, respectively, to the identities
\begin{align}
	h_1 \rhd \Bigl( \bigl( h^* \lhu \tau(h_2) \bigr) \rhd v \Bigr) & = \bigl( \sigma(h_1) \rhu h^* \bigr) \rhd (h_2 \rhd v), \label{eq:yd-dual-lr} \\
	\Bigl( \bigl( \sigma(h_2) \rhu h^* \bigr) \rhd v \Bigr) \lhd h_1 & = \bigl( h^* \lhu \tau(h_1) \bigr) \rhd (v \lhd h_2), \label{eq:yd-dual-rr} \\
	\Bigl( v \lhd \bigl( \sigma(h_1) \rhu h^* \bigr) \Bigr) \lhd h_2 & = (v \lhd h_1) \lhd \bigl( h^* \lhu \tau(h_2) \bigr), \label{eq:yd-dual-rl} \\
	h_2 \rhd \Bigl( v \lhd \bigl( h^* \lhu \tau(h_1) \bigr) \Bigr) & = (h_1 \rhd v) \lhd (\sigma(h_2) \rhu h^*). \label{eq:yd-dual-ll}
\end{align}
We define $\widehat{\YDlr}(\sigma,\tau)$ to be the category whose objects are $V$ equipped with a left $H$-module structure and a left $H^*$-module structure satisfying \eqref{eq:yd-dual-lr}, and whose morphisms are both $H$-linear and $H^*$-linear.
If $H$ is a Hopf algebra with bijective antipode, then $\widehat{\YDlr}(\sigma, \tau)$ is equivalent to the category of left modules over the algebra $H^* \#_{\sigma,\tau} H$ of \Cref{sssec:twisted-Drinfeld-double}.

If $H$ is finite-dimensional, these identities are precisely the conditions for $V$ to be a module over the twisted Drinfeld double $D_{\sigma, \tau}(H)$ of \cite[Section~2]{PS07}; 
for $\sigma = \tau = \id_H$, this recovers the classical description of Yetter--Drinfeld modules as modules over the Drinfeld double $D(H)$.

There is a faithful exact functor $F \colon \YDlr(\sigma,\tau) \to \widehat{\YDlr}(\sigma,\tau)$. An object $V \in \widehat{\YDlr}(\sigma,\tau)$ is isomorphic to $F(U)$ for some $U \in \YDlr(\sigma,\tau)$ if and only if $V$ is rational as an $H^*$-module; that is, its module structure is induced by a right $H$-comodule structure via $h^* \rhd v = v_0 \langle h^*, v_{1} \rangle$.
The categories $\widehat{\YDrr}(\sigma,\tau)$, $\widehat{\YDrl}(\sigma,\tau)$, and $\widehat{\YDll}(\sigma,\tau)$ can be defined analogously.


The following duality is \cite[Proposition~3.7]{PS07} when $\sigma$ and $\tau$ are Hopf algebra automorphisms.
The same formulas remain valid for the bialgebra endomorphisms used here.

\begin{lem}\label{lem:dual-YD}
	Suppose that $V \in \YDrl(\sigma, \tau)$ is finite-dimensional.
	Then $V^* = \Hom(V, \Bbbk)$ is an object of $\YDlr(\sigma, \tau)$ with
	\[ (h \rhd f)(v) = f(v \lhd h), \quad \text{and} \quad f_{0}(v) \otimes f_1 = f(v_{0}) \otimes v_{-1}. \]
\end{lem}
\begin{proof}
	By the definition of the coaction on $V^*$ and the compatibility condition \eqref{eq:yd-rl} for $V$, we have
		\begin{align*}
			(h_1 \rhd f_0)(v) \otimes \tau(h_2)f_{1} & = f_0(v\lhd h_1) \otimes \tau(h_2)f_{1} \\
			& = f\bigl( (v\lhd h_1)_0 \bigr) \otimes \tau(h_2)(v\lhd h_1)_{-1} \\
			& \xlongequal{\eqref{eq:yd-rl}} f(v_{0}\lhd h_2) \otimes v_{-1}\sigma(h_1) \\
			& = (h_2\rhd f)(v_{0}) \otimes v_{-1}\sigma(h_1) \\
			& = (h_2\rhd f)_{0}(v) \otimes (h_2\rhd f)_{1}\sigma(h_1).
		\end{align*}
	This is precisely the compatibility condition \eqref{eq:yd-lr}, so $V^* \in \YDlr(\sigma, \tau)$.
\end{proof}

\Cref{lem:dual-YD} relates the two twisted Yetter--Drinfeld module categories with the same Hopf endomorphisms $\sigma$ and $\tau$ via duality.


\subsubsection{The twisted Drinfeld double}\label{sssec:twisted-Drinfeld-double}

Let $H$ be a Hopf algebra with bijective antipode.
Then the abelian category $\widehat{\YDlr}(\sigma,\tau)$ is equivalent to the category of modules over an algebra known as the twisted Drinfeld double.

Recall that, given a bialgebra $H$, a left $H$-module algebra $A$, and a left $H$-comodule algebra $B$, the Takeuchi smash product $A \# B$, 
is the associative algebra with underlying $\Bbbk$-space $A \otimes B$, multiplication
\[ (a \# b)(a' \# b') := a (b_{-1} \rhd a') \# b_{0} b', \quad a, a' \in A,\ b, b' \in B, \]
and unit $1_A \# 1_B$.

Now let $H$ be a Hopf algebra with bijective antipode $S$.
Then $H$ is a left $H \otimes H^{\op}$-comodule algebra via
\[ \rho \colon H \To (H \otimes H^{\op}) \otimes H, \quad h \mapsto \bigl( h_1 \otimes S^{-1}(h_3) \bigr) \otimes h_2. \]
Moreover, for any two Hopf algebra endomorphisms $\sigma, \tau \colon H \to H$, the dual algebra $H^*$ is a left $H \otimes H^{\op}$-module algebra 
via
\[ (h \otimes l) \rhu f = \sigma(h) \rhu f \lhu \tau(l). \]
We denote by $H^* \#_{\sigma, \tau} H$ the Takeuchi smash product associated with the left $H \otimes H^{\op}$-module algebra $H^*$ and the left $H \otimes H^{\op}$-comodule algebra $H$ defined above.
When $H$ is finite-dimensional, \cite[Section~2]{PS07} proves that $\widehat{\YDlr}(\sigma,\tau)$ (respectively, $\widehat{\YDrl}(\sigma,\tau)$) is equivalent to the category of left (respectively, right) $H^* \#_{\sigma,\tau} H$-modules. The same equivalence holds for the algebra constructed above, where $\sigma$ and $\tau$ are Hopf algebra endomorphisms and $H$ need not be finite-dimensional.

When $\sigma=\tau=\id_H$ and $H$ is finite-dimensional, $H^*\# H$ is the underlying algebra of the Drinfeld double $D(H)=H^{*\cop}\bowtie H$ \cite{Dri87}.
Thus $H^*\#_{\sigma,\tau}H$ is a twisted form of that algebra.

For a general Hopf algebra, $H^*$ need not be a coalgebra. The Sweedler finite dual
\[
H^\circ
=
\{\,f\in H^*\mid f(I)=0
\text{ for some finite-codimensional ideal }I\subset H\,\}
\]
is a Hopf algebra and is stable under the untwisted $H\otimes H^{\op}$-action above. Passing to the coopposite coalgebra does not change the multiplication, so the Drinfeld double $D(H)=(H^\circ)^{\cop}\bowtie H$ of \cite{DT94, Maj95} satisfies
\[
D(H)\cong H^\circ\# H\hookrightarrow H^*\# H.
\]
If $H$ is finite-dimensional, then $H^\circ=H^*$ and the inclusion is an equality.

\subsubsection{The Drinfeld center of a bimodule category}

The category of twisted Yetter--Drinfeld modules is the Drinfeld center of a suitable bimodule category. The $S^4$ formula for finite tensor categories follows from the fact that the distinguished invertible object lies in the Drinfeld center of the bimodule category whose left action is pulled back along the quadruple dual functor.

Let $\cC$ be a monoidal category with associativity constraint
$a_{X, Y, Z} \colon (X \otimes Y) \otimes Z \stackrel{\cong}{\to} X \otimes (Y \otimes Z)$,
and let $\cM$ be a $\cC$-bimodule category.
Recall that the \emph{Drinfeld center} $\cZ_{\cC}(\cM)$ of $\cM$ \cite[Definition~2.1]{GNN09} is the category whose objects are pairs $(M,\gamma_M)$, where $M\in\cM$ and
\[ \gamma_M = \{\, \gamma_{M,X} \colon X \otimes M \stackrel{\cong}{\to} M \otimes X \,\}_{X \in \cC} \]
is a family of isomorphisms that is natural in $X$ and satisfies the half-braiding (hexagon) condition.
A morphism from $(M, \gamma_M)$ to $(M', \gamma_{M'})$ is a morphism $f \in \Hom_{\cM}(M, M')$
such that, for every $X \in \cC$,
\[
(f \otimes \id_{X}) \circ \gamma_{M,X} = \gamma_{M',X} \circ (\id_X \otimes f).
\]
When $\cM = \cC$, this construction coincides with the ordinary Drinfeld center $\cZ(\cC)$.

Let $(\cM^H, \otimes, \Bbbk^{\eta})$ denote the monoidal category of right $H$-comodules, and let $\sigma$ and $\tau$ be Hopf algebra endomorphisms.
Let $\cM^H(\sigma, \tau)$ denote the $\cM^H$-bimodule category $\cM^H$ whose left action $\otimes^{\ell}$ is twisted by $\tau$ and whose right action $\otimes^{r}$ is twisted by $\sigma$; that is,
\[ N \otimes^{\ell} M := N^{\tau} \otimes M, \quad M \otimes^{r} N := M \otimes N^{\sigma}, \quad N \in \cM^H, \ M \in \cM^H. \]
The Drinfeld center $\cZ_{\cM^H}\bigl(\cM^H(\sigma,\tau)\bigr)$ of the $\cM^H$-bimodule category $\cM^H(\sigma,\tau)$ is equivalent to the category of twisted Yetter--Drinfeld modules $\YDlr(\sigma,\tau)$.

\begin{rem}
	Let $\cC$ be a finite tensor category, and let $X_{\rho} \in \cC$ be its distinguished invertible object. Etingof, Nikshych, and Ostrik \cite{ENO04} proved that, for every $V \in \cC$, there is an isomorphism $\gamma_V \colon V^{****} \otimes X_{\rho} \xrightarrow{\;\cong\;} X_{\rho} \otimes V$.
	Moreover, $(X_{\rho}, \gamma = \{\gamma_V\}_{V \in \cC})$ defines an object of the Drinfeld center $\cZ(\cM)$, where $\cM$ is the $\cC$-bimodule category whose underlying category is $\cC$, whose right action is the usual one, and whose left action is pulled back along the quadruple dual functor; see \cite[Remark~4.13]{FSS20}.
	
	If $\cC$ is the category of finite-dimensional comodules over a finite-dimensional Hopf algebra $H$, then $\cM$ is equivalent to the full subcategory of $\cM^H(\id_H, S^4)$ consisting of finite-dimensional objects, and $X_{\rho} \cong \rcoint$ as right $H$-comodules.
	However, when $H$ is infinite-dimensional, this conclusion cannot be applied directly, as the category of finite-dimensional $H$-comodules need not be a finite tensor category.
\end{rem}

\section{\texorpdfstring{$S^4$-formula for AS (co-)Gorenstein Hopf algebras}{S4 formula for AS (co-)Gorenstein Hopf algebras}}\label{sec:S4}

We prove the $S^4$ formula for AS-Gorenstein Hopf algebras, and then the corresponding formula for AS-co-Gorenstein Hopf algebras.

\subsection{\texorpdfstring{$S^4$-formula for AS-Gorenstein Hopf algebras}{S4 formula for AS-Gorenstein Hopf algebras}}

In this subsection, we derive a formula for the fourth power of the antipode of an AS-Gorenstein Hopf algebra by constructing an $(\id_H, S^4)$-twisted right-left Yetter--Drinfeld module structure on its space of homological integrals.
Our strategy is to construct a resolution of the trivial module ${}_{\vep}\Bbbk$ by free $H$-modules in a suitable category of Yetter--Drinfeld modules, so that the cohomology $\Ext_H^*({}_{\vep}\Bbbk, H)$ inherits a canonical twisted Yetter--Drinfeld structure.
We first construct a functor $H \boxtimes -$ from the category of left $H$-comodules to the category of twisted Yetter--Drinfeld modules and show that it is left adjoint to the forgetful functor (\Cref{lem:left-adjoint}). This adjunction yields the desired $\cF$-resolutions (\Cref{prop:resolution}).
We then construct compatible $H$-module and $H^*$-module structures on the relevant Ext groups that are independent of the choice of $\cF$-resolution (\Cref{prop:Ext-structure}).
Under the rationality hypotheses, these structures yield twisted Yetter--Drinfeld modules.
Combining these preparatory results, we prove the main result of this subsection (\Cref{thm:main}) and conclude by comparing our formula with the classical distinguished group-like element in the finite-dimensional case.

Let $H$ be a Hopf algebra with antipode $S$, and let $\sigma$ and $\tau$ be Hopf algebra endomorphisms of $H$.
For any left $H$-comodule $N$, multiplication on the first factor endows $H \otimes N$ with a left $H$-module structure. Moreover, the linear map
\[ H \otimes N \To H \otimes H \otimes N, \quad h \otimes n \mapsto \tau(h_1)n_{-1}\sigma(S(h_3)) \otimes h_2 \otimes n_0 \]
endows $H \otimes N$ with a left $H$-comodule structure.
A direct verification shows that the module and comodule structures on $H \otimes N$ satisfy the compatibility condition \eqref{eq:yd-ll}.
Hence $H \otimes N$ is a $(\sigma, \tau)$-twisted left-left Yetter--Drinfeld module, which we denote by $H \boxtimes N$.
This construction produces a functor
\[ H \boxtimes - \colon {^H\!\!\cM} \To {\YDll}(\sigma, \tau), \quad N \mapsto H \boxtimes N. \]

This functor admits a right adjoint that forgets the module
structure. The following lemma is a twisted version of \cite[Corollary 2.8]{CMZ97}. As usual, we omit the forgetful functor from the notation.
\begin{lem}\label{lem:left-adjoint}
	The functor $H \boxtimes - \colon {^H\!\!\cM} \to {\YDll(\sigma, \tau)}$ is left adjoint to the forgetful functor ${\YDll}(\sigma, \tau) \to {^H\!\cM}$; that is, there is a natural isomorphism
	\[ \Phi_{N, V} \colon \Hom_{^H\!\!\cM}(N, V) \To \Hom_{\YDll(\sigma, \tau)}(H \boxtimes N, V), \quad f \mapsto \big(h \otimes n \mapsto h \rhd f(n) \big). \]
\end{lem}

Denote by $\cF$ the full subcategory of $\YDll(\sigma,\tau)$ consisting of objects of the form $H \boxtimes N$, where $N$ is a left $H$-comodule.

\begin{prop}\label{prop:resolution}
	Let $V \in {\YDll}(\sigma, \tau)$.
	\begin{enumerate}[label=\textup{(\arabic*)},leftmargin=2.2em]
		\item The morphism $\varphi_V := \Phi_{V,V}(\id_V)$ is an \emph{$\cF$-precover} of $V$; that is, for every 
        $F \in \cF$, the map
		\[ (\varphi_V)_* \colon \Hom_{{\YDll}(\sigma,\tau)}(F, H \boxtimes V) \To \Hom_{{\YDll}(\sigma,\tau)}(F, V) \]
		is surjective.
		\item There exists an \emph{$\cF$-resolution} of $V$; that is, a (not necessarily exact) complex in ${\YDll}(\sigma,\tau)$,
        \[\cdots \To F_1 \To F_0 \To V \To 0\]
        such that each $F_i \in \cF$ and the complex becomes exact after applying $\Hom_{{\YDll}(\sigma,\tau)}(F,-)$ for any $F \in \cF$.
		\item In particular, the $\cF$-resolution of $V$ may be chosen to be an exact sequence
		\[ \cdots \stackrel{d_{i+1}}{\To} H \boxtimes N_i \stackrel{d_i}{\To} \cdots \stackrel{d_2}{\To} H \boxtimes N_1 \stackrel{d_1}{\To} H \boxtimes N_0 \stackrel{d_0}{\To} V \To 0 \]
		in ${\YDll}(\sigma,\tau)$, where, for each $i \geqslant 0$, $N_i$ is an $H$-comodule and the restriction $d_i|_{1 \boxtimes N_i}$ is injective.
		\item If $V$ is of type $\mathrm{FP}_{n} (n \geqslant 0)$ as a left $H$-module, then there exists an exact sequence
		\[ \cdots \To H \boxtimes N_i \To \cdots \To H \boxtimes N_1 \To H \boxtimes N_0 \To V \To 0 \]
		in ${\YDll}(\sigma,\tau)$, where $N_i$ is finite-dimensional for $0 \leqslant i \leqslant n$.
		This sequence is a free resolution of the underlying $H$-module $V$, although it is not claimed to be an $\cF$-resolution in the sense of (2).
	\end{enumerate}
\end{prop}
\begin{proof}
	(1)
	Since $F \in \cF$, there exists an $H$-comodule $N$ such that $F = H \boxtimes N$.
	For any $f \in \Hom_{{\YDll}(\sigma,\tau)}(F, V)$, define $\widetilde{f}:  F \to H \boxtimes V, \, h \boxtimes n \mapsto h \boxtimes f(1 \boxtimes n)$. Then $\widetilde{f} \in \Hom_{{\YDll}(\sigma,\tau)}(F, H \boxtimes V)$.
	For any $h \in H$ and $n \in N$,
	\[ (\varphi_V \widetilde{f})(h \boxtimes n) = \varphi_V\bigl( h \boxtimes f(1 \boxtimes n) \bigr) = h \rhd f(1 \boxtimes n) = f(h \boxtimes n), \]
	where the last equality follows from the $H$-linearity of $f$.
   Hence $\varphi_V \widetilde{f} = f$, and $(\varphi_V)_*$ is surjective.
	
	(2) This follows directly from (1) and \cite[Proposition~8.1.3]{EJ11}.
	
	(3)
	Set $N_0 = V$ and $d_0 = \varphi_{N_0} = \Phi_{V,V}(\id_V)$.
	Then $d_0$ is surjective, and its restriction $d_0|_{1 \boxtimes N_0}$ is injective.
	Applying the same construction to $N_1=\Ker(d_0)$ and proceeding inductively, we obtain the desired exact sequence.
	
	(4)
	Since $V$ is of type $\mathrm{FP}_{n}$ as a left $H$-module, ${}_HV$ is finitely generated.
	Fix a set of generators $\{x_1, \dots, x_m\}$ of ${}_HV$.
	By the Finiteness Theorem \cite[5.1.1]{Mon93}, there exists a finite-dimensional $H$-subcomodule $N_0$ of $V$ containing $\{x_1, \dots, x_m\}$.
	Let $\iota \colon N_0 \hookrightarrow V$ denote the inclusion.
	Consequently, there is a morphism
	$d_0 := \Phi_{N_0, V}(\iota) \colon H \boxtimes N_0 \to V$ in
	$\YDll(\sigma,\tau)$.
	Since $H N_0 = V$, the map $d_0$ is surjective.
	Suppose now that $n \geqslant 1$.
	The module $H \boxtimes N_0$ is a finitely generated free left $H$-module and is therefore of type $\mathrm{FP}_{n-1}$ as a left $H$-module.
	By Schanuel's lemma \cite[Proposition~3.12]{Rot09}, $\mathrm{Ker}(d_0)$ is also of type $\mathrm{FP}_{n-1}$ as a left $H$-module.
	The desired conclusion follows by induction on $n$.
\end{proof}

Bichon \cite{Bic13} used the untwisted case of the resolution in \Cref{prop:resolution} (3) to compute bialgebra cohomology in the sense of \cite{GS90}.

Since each $H \boxtimes N_i$ is free as a left $H$-module, such a resolution can be used to compute $\Ext$-groups.
To endow these $\Ext$-groups with a twisted Yetter--Drinfeld module structure, we first study the Hom spaces $\Hom_H(V, X)$ in a general setting.

\begin{lem}\label{lem:Hom-structure}
  Let $V \in \YDll(\sigma, \tau)$, and let $X\in {_H^H\!\cM^H_H}$ be a Hopf bimodule \cite[Section~3]{Sch94}.
	\begin{enumerate}[label=\textup{(\arabic*)},leftmargin=2.2em]
		\item The Hom space $\Hom_H(V, X)$ carries a natural right $H^*$-module structure given by
		\begin{equation}\label{eq:Hom-H^*-action}
			(f \lhd h^*)(v) = f(v_0)_0 \bigl\langle h^*, \sigma\Bigl( S^3\bigl( f(v_0)_{1} \bigr) \Bigr) S(v_{-1}) \tau\bigl( f(v_0)_{-1} \bigr) \bigr\rangle. 
		\end{equation}
		\item $\Hom_H(V, X)$ also carries a natural right $H$-module structure induced by the right $H$-module structure on $X$, namely,
		\begin{equation}\label{eq:Hom-H-action}
			(f \lhd h)(v) = f(v) \lhd h.
		\end{equation}
		This right $H$-module structure and the right $H^*$-module structure on $\Hom_H(V, X)$ given in (1) satisfy \eqref{eq:yd-dual-rl} for the pair $(\tau, \sigma S^4)$ of Hopf algebra endomorphisms.
		\item If $\Hom_H(V, X)$, equipped with the structure defined by \eqref{eq:Hom-H^*-action}, is a rational right $H^*$-module, then $\Hom_H(V, X)$, endowed with the right $H$-module structure defined by \eqref{eq:Hom-H-action}, is an object of $\YDrl(\tau, \sigma S^4)$.
		\item If $V$ is finitely generated as a left $H$-module, then $\Hom_H(V, X)$ is a rational right $H^*$-module and hence, by (3), lies in $\YDrl(\tau, \sigma S^4)$.
	\end{enumerate}
\end{lem}
\begin{proof}
	(1) It is easy to see that equation \eqref{eq:Hom-H^*-action} induces a right $H^*$-module structure on $\Hom(V, X)$.
    Using the antipode, the compatibility condition \eqref{eq:yd-ll} for a $(\sigma,\tau)$-twisted left-left Yetter--Drinfeld module can equivalently be written as
    \begin{equation}\label{eq:yd-ll-S}
    (h\rhd v)_{-1}\otimes(h\rhd v)_0
    =
    \tau(h_1)v_{-1}\sigma\bigl(S(h_3)\bigr)
    \otimes h_2\rhd v_0
    \quad(h\in H,\ v\in V).
    \end{equation}
	To prove (1), it suffices to show that $f \lhd h^*$ is $H$-linear for any $f \in \Hom_H(V, X)$ and $h^* \in H^*$.
	This is verified as follows:
	\begin{align*}
		& (f \lhd h^*)(h \rhd v) \\
		& \xlongequal{\eqref{eq:Hom-H^*-action}} f\bigl( (h \rhd v)_0 \bigr)_0 \bigl\langle h^*, \sigma\Bigl( S^3\bigl( f( (h \rhd v)_0 )_{1} \bigr) \Bigr) S\bigl( (h \rhd v)_{-1} \bigr) \tau\bigl( f( (h \rhd v)_0 )_{-1} \bigr)  \bigr\rangle \\
		& \xlongequal{\eqref{eq:yd-ll-S}} f(h_2 \rhd v_0)_0 \bigl\langle h^*, \sigma\Bigl( S^3\bigl( f(h_2 \rhd v_0)_{1} \bigr) \Bigr) S\Bigl( \tau(h_1) v_{-1} \sigma\bigl( S(h_3) \bigr) \Bigr) \tau\bigl( f(h_2 \rhd v_0)_{-1} \bigr) \bigr\rangle \\
		& = h_3 \rhd f(v_0)_0 \bigl\langle h^*, \sigma\Bigl( S^3 \bigl( f(v_0)_{1} \bigr) S^3(h_4) S^2(h_5) \Bigr) S(v_{-1}) \tau\bigl( S(h_1) h_2 f(v_0)_{-1} \bigr) \bigr\rangle \\
		& = h \rhd f(v_0)_0 \bigl\langle h^*, \sigma\Bigl( S^3\bigl( f(v_0)_{1} \bigr) \Bigr) S(v_{-1}) \tau\bigl( f(v_0)_{-1} \bigr) \bigr\rangle \\
		& = h \rhd (f \lhd h^*)(v).
	\end{align*}
	
	(2)
	Obviously, the formula \eqref{eq:Hom-H-action} defines a right $H$-module structure on $\Hom_H(V, X)$.
	A direct computation gives
	\begin{align*}
		& \Bigl( \bigl( f \lhd (\tau(h_1) \rhu h^*) \bigr) \lhd h_2 \Bigr) (v) \\
		& = \bigl( f \lhd (\tau(h_1) \rhu h^*) \bigr) (v) \lhd h_2 \\
		& \xlongequal{\eqref{eq:Hom-H^*-action}} f(v_0)_0 \lhd h_2 \bigl\langle \tau(h_1) \rhu h^*, \sigma\Bigl( S^3\bigl( f(v_0)_{1} \bigr) \Bigr) S(v_{-1}) \tau\bigl( f(v_0)_{-1} \bigr) \bigr\rangle \\
		& = f(v_0)_0 \lhd h_2 \bigl\langle h^*, \sigma\Bigl(S^4(h_4) S^3(h_3) S^3\bigl( f(v_0)_{1} \bigr) \Bigr) S(v_{-1}) \tau\bigl( f(v_0)_{-1} \bigr) \tau(h_1) \bigr\rangle \\
		& = (f \lhd h_1)(v_0)_0 \bigl\langle h^* \lhu \sigma\bigl( S^4(h_2) \bigr), \sigma\Bigl( S^3\bigl( (f \lhd h_1)(v_0)_{1} \bigr) \Bigr) S(v_{-1}) \tau\bigl( (f \lhd h_1)(v_0)_{-1} \bigr) \bigr\rangle \\
		& = \Bigl( (f \lhd h_1) \lhd \bigl( h^* \lhu \sigma(S^4(h_2)) \bigr) \Bigr) (v).
	\end{align*}
	Therefore, \eqref{eq:yd-dual-rl} holds for the pair $(\tau, \sigma S^4)$ of Hopf algebra endomorphisms.
	
	(3)
	It follows from (2) and the equivalence between \eqref{eq:yd-rl} and \eqref{eq:yd-dual-rl} when the right $H^*$-module structure is induced by a left $H$-comodule structure.
	
	(4)
	By \Cref{prop:resolution}(4), there exist a finite-dimensional $H$-comodule $N$ and a surjective morphism $\pi: H \boxtimes N \to V$.
	Hence $\Hom_H(V, X)$ is an $H^*$-submodule of $\Hom_H(H \boxtimes N, X)$.
	Thus, it suffices to show that $\Hom_H(H \boxtimes N, X)$ is a rational right $H^*$-module.
	Let $N^*$ be the $\Bbbk$-linear dual $\Hom_{\Bbbk}(N, \Bbbk)$, endowed with the left $H$-comodule structure determined by $f_{-1} \langle f_0, n \rangle = S(n_{-1}) \langle f, n_0 \rangle$ for any $f \in N^*$ and $n \in N$.
	Consider the $\Bbbk$-linear map
	\[ \Theta: N^* \otimes X \To \Hom_H(H \boxtimes N, X), \quad f \otimes x \mapsto \bigl( h \boxtimes n \mapsto f(n) h \rhd x \bigr), \]
	where $N^* \otimes X$ is endowed with the right $H^*$-module structure induced by the left $H$-comodule structure
	\[ f \otimes x \mapsto \sigma (S^3(x_1)) f_{-1} \tau(x_{-1}) \otimes f_0 \otimes x_0. \]
	The map $\Theta$ is bijective.
	For $g = \Theta(f \otimes x)$, we obtain
	\begin{align*}
		& \bigl( g \lhd h^* \bigr) (h \boxtimes n) \\
		& = g (h_2 \boxtimes n_0)_0 \langle h^*, \sigma \Bigl( S^3 \bigl( g (h_2 \boxtimes n_0)_1 \bigr) \Bigr) S\bigl( \tau(h_1) n_{-1} \sigma(S(h_3)) \bigr) \tau\bigl( g (h_2 \boxtimes n_0)_{-1} \bigr) \rangle \\
		& = f(n_0) (h_2 \rhd x)_0 \langle h^*, \sigma \Bigl( S^3 \bigl( (h_2 \rhd x)_1 \bigr) \Bigr) S\bigl( \tau(h_1) n_{-1} \sigma(S(h_3)) \bigr) \tau\bigl( (h_2 \rhd x)_{-1} \bigr) \rangle \\
		& = f(n_0) (h_3 \rhd x_0) \langle h^*, \sigma \bigl( S^3 (h_4x_1) \bigr) S\bigl( \tau(h_1) n_{-1} \sigma(S(h_5)) \bigr) \tau(h_2x_{-1}) \rangle \\
		& = f(n_0) (h \rhd x_0) \langle h^*, \sigma \bigl( S^3 (x_1) \bigr) S(n_{-1}) \tau(x_{-1}) \rangle \\
		& = f_0(n) (h \rhd x_0) \langle h^*, \sigma \bigl( S^3 (x_1) \bigr) f_{-1} \tau(x_{-1}) \rangle \\
		& = \Theta(f_0 \otimes x_0)(h \boxtimes n) \langle h^*, \sigma \bigl( S^3 (x_1) \bigr) f_{-1} \tau(x_{-1}) \rangle \\
		& = \Theta\bigl( (f \otimes x) \lhd h^* \bigr) (h \boxtimes n).
	\end{align*}
	Therefore, $\Theta$ is $H^*$-linear and hence an isomorphism of right $H^*$-modules.
\end{proof}

We use the following comparison lemma, whose proof is routine and is therefore omitted (cf. \cite[Section~8.1]{EJ11}).

\begin{lem}\label{lem:comparison}
	Let $U, V \in \YDll(\sigma, \tau)$, let $\cdots \to P_1 \to P_0 \to U \to 0$ be a complex in $\YDll(\sigma, \tau)$ with $P_i \in \cF$ for all $i \geqslant 0$, and let $F_{\bullet}$ be an $\cF$-resolution of $V$.
	Then every morphism $f \colon U \to V$ in $\YDll(\sigma, \tau)$ lifts to a chain map $P_{\bullet} \to F_{\bullet}$, unique up to homotopy.
\end{lem}

By \Cref{lem:comparison}, every morphism in $\YDll(\sigma, \tau)$ lifts to a chain map of $\cF$-resolutions, and the lift is unique up to homotopy. The following proposition uses this lift.

\begin{prop}\label{prop:Ext-structure}
	Let $V \in \YDll(\sigma, \tau)$ and let $X \in {_H^H\!\cM^H_H}$ be a Hopf bimodule.
	For every $i \geqslant 0$, the following statements hold.
	\begin{enumerate}[label=\textup{(\arabic*)},leftmargin=2.2em]
		\item The space $\Ext_H^i(V, X)$ carries a natural right $H^*$-module structure induced by \eqref{eq:Hom-H^*-action} and a natural right $H$-module structure induced by \eqref{eq:Hom-H-action}.
		Moreover, both structures are independent of the choice of the $\cF$-resolution of $V$ constructed in \Cref{prop:resolution}(3).
		\item The space $\Ext_H^i(V, X)$ is an object of $\widehat{\YDrl}(\tau, \sigma S^4)$.
		Furthermore, if it is rational as a right $H^*$-module, then it is an object of $\YDrl(\tau, \sigma S^4)$.
	\end{enumerate}
\end{prop}
\begin{proof}
	(1)
	Let $F_{\bullet}$ be an $\cF$-resolution of $V$ constructed in \Cref{prop:resolution}(3).
	Since $F_{\bullet}$ is also a projective resolution of $V$ as an $H$-module, the cohomology groups $\mathrm{H}^i\bigl(\Hom_H(F_{\bullet}, X)\bigr)$ coincide, as $\Bbbk$-vector spaces, with the usual Ext groups $\Ext_H^i(V, X)$ in the category of $H$-modules.
	By \Cref{lem:Hom-structure}, the complex $\Hom_H(F_{\bullet}, X)$ is both an $H$-complex and an $H^*$-complex.
	Hence, its cohomology groups inherit right $H$-module and right $H^*$-module structures.
	
	Suppose that $F'_{\bullet}$ is another $\cF$-resolution of $V$.
	By \Cref{lem:comparison}, there exists a chain map $f_{\bullet} \colon F'_{\bullet} \to F_{\bullet}$ that lifts $\id_V$ and is a homotopy equivalence. 
	Thus
	\[ \Hom_H(f_{\bullet}, X) \colon \Hom_H(F_{\bullet}, X) \To \Hom_H(F'_{\bullet}, X) \]
	is also a homotopy equivalence.
	Since each $\Hom_H(f_i, X)$ is both $H$-linear and $H^*$-linear, the two induced module structures on the cohomology groups agree.

  (2)
  By \Cref{lem:Hom-structure}(2), each term of $\Hom_H(F_{\bullet}, X)$ satisfies \eqref{eq:yd-dual-rl} for the pair $(\tau, \sigma S^4)$. The differentials are morphisms of both module structures, so the identity passes to the cohomology, and $\Ext_H^i(V, X)$ lies in $\widehat{\YDrl}(\tau, \sigma S^4)$. If this right $H^*$-module is rational, then $\Ext_H^i(V, X)$ lies in $\YDrl(\tau, \sigma S^4)$, as in \Cref{lem:Hom-structure}(3).
\end{proof}

\begin{lem}\label{lem:FP-rational}
	Let $H$ be a Hopf algebra with antipode $S$, let $V \in \YDll(\sigma, \tau)$, and let $X \in {_H^H\!\cM^H_H}$ be a Hopf bimodule.
	If $V$ is of type $\mathrm{FP}_{d}$ as an $H$-module, then $\Ext^i_H(V, X)$ is a rational right $H^*$-module for every $i \leqslant d$.
\end{lem}
\begin{proof}
	By \Cref{prop:resolution}(4), there exists an exact sequence $\cdots \to P_1 \to P_0 \to V \to 0$ in $\YDll(\sigma, \tau)$, where $P_i = H \boxtimes N_i$ and $N_i$ is finite-dimensional for $i \leqslant d$.
	By \Cref{lem:Hom-structure}(4), $\Hom_H(P_i, X)$ is rational as a right $H^*$-module for $i \leqslant d$.
	Since rational modules are closed under submodules and quotients, $\mathrm{H}^i\bigl(\Hom_H(P_{\bullet}, X)\bigr)$ is rational for $i \leqslant d$.
	
	Let $F_{\bullet}$ be an $\cF$-resolution of $V$ constructed in \Cref{prop:resolution}(3).
	By \Cref{lem:comparison}, there exists a chain map $f_{\bullet} \colon P_{\bullet} \to F_{\bullet}$ lifting $\id_V$ in ${\YDll}(\sigma,\tau)$.
	Since both $P_{\bullet}$ and $F_{\bullet}$ are projective resolutions of $V$ as left $H$-modules, $f_{\bullet}$ is a homotopy equivalence of complexes of $H$-modules.
	Moreover, each $\Hom_H(f_i, X)$ is both $H$-linear and $H^*$-linear, so
	$f_{\bullet}$ induces isomorphisms
	\[ \mathrm{H}^i\bigl(\Hom_H(P_{\bullet}, X)\bigr) \cong \mathrm{H}^i\bigl(\Hom_H(F_{\bullet}, X)\bigr) = \Ext_H^i(V, X) \]
	of right $H$-modules and right $H^*$-modules.
	In particular, $\Ext_H^i(V, X)$ is rational as a right $H^*$-module for $i \leqslant d$.
\end{proof}


\begin{lem}\cite[Proposition~3.5.3]{HR74}\label{lem:coreflexive-rational}
  Let $C$ be a coalgebra, and let $C_0$ denote its coradical. Then $C_0$ is coreflexive if and only if every finite-dimensional simple right $C^*$-module is rational, and if and only if every finite-dimensional simple left $C^*$-module is rational. 
\end{lem}

\begin{lem}\cite[Corollary~3.6 and~3.12]{Rad73}\label{lem:coreflexive}
	Let $C$ be a coalgebra, and let $\cS(C)$ denote the set of simple subcoalgebras of $C$.
	Suppose that $\Bbbk$ is infinite.
	If $|\cS(C)| \leqslant |\Bbbk|$, then the coradical $C_0$ is coreflexive.
\end{lem}
Further conditions ensuring that the coradical is coreflexive are given in \cite[Theorem~3.7.5]{HR74}.

We are now in a position to state and prove the main result: the space of left homological integrals of a left AS-Gorenstein Hopf algebra carries a canonical twisted Yetter--Drinfeld structure, from which the $S^4$-formula follows.
\begin{thm}\label{thm:main}
	Let $H$ be a left AS-Gorenstein Hopf algebra of dimension $d$ with antipode $S$.
	Suppose that either the coradical of $H$ is coreflexive or ${}_{\vep}\Bbbk$ is of type $\mathrm{FP}_{d}$.
	\begin{enumerate}[label=\textup{(\arabic*)},leftmargin=2.2em]
		\item The space of left homological integrals  $\lint$ is an $(\id_H, S^4)$-twisted right-left Yetter--Drinfeld module.
		Thus, there exist a group-like element $g$ and a character $\alpha \colon H \to \Bbbk$ such that
		\[
		\lint \cong \Bbbk g \ \text{as a left $H$-comodule}; \,
		\lint \cong \Bbbk_\alpha \ \text{as a right $H$-module}.
		\]
		\item For all $h\in H$,
		\begin{equation}\label{eq:main-formula}
			S^4(h) = g (\alpha \rhu h \lhu \alpha^{-1}) g^{-1}.
		\end{equation}
	\end{enumerate}
\end{thm}
\begin{proof}
	(1)
	Note that $\Bbbk$ is a left-left Yetter--Drinfeld module with the trivial left $H$-module structure and the trivial left $H$-comodule structure.
	Applying \Cref{prop:Ext-structure} with $\sigma = \tau = \id_H$, we obtain
	\[ \Ext^i_H(\Bbbk, H) \in \widehat{\YDrl}(\id_H, S^4) \, \text{ for all } i \geqslant 0. \]
	Since $H$ is left AS-Gorenstein of dimension $d$, the space $\lint = \Ext^d_H(\Bbbk, H)$ is one-dimensional and therefore simple.
	If the coradical of $H$ is coreflexive, then $\lint$ is rational as a right $H^*$-module by \Cref{lem:coreflexive-rational}.
	If ${}_{\vep}\Bbbk$ is of type $\mathrm{FP}_{d}$, then $\lint$ is also rational as a right $H^*$-module by \Cref{lem:FP-rational}.
	Thus, in either case, $\lint$ is rational as a right $H^*$-module and hence $\lint \in \YDrl(\id_H, S^4)$.
	Consequently, as a left $H$-comodule, $\lint$ is of the form $\Bbbk g$ for some group-like element $g \in G(H)$ and, as a right $H$-module, of the form $\Bbbk_{\alpha}$ for some character $\alpha \in G(H^{\circ})$.
	
	(2)
	Fix a nonzero element $t \in \Ext^d_H(\Bbbk, H)$.
	By the compatibility condition \eqref{eq:yd-rl} with $(\sigma, \tau) = (\id_H, S^4)$, we obtain
	\[ g h_1 \otimes t \alpha(h_2) = S^4(h_2) (t \lhd h_1)_{-1} \otimes (t \lhd h_1)_0 = S^4(h_2) t_{-1} \otimes t_0 \alpha(h_1) = S^{4}(h_2) g \otimes t \alpha(h_1). \]
	Since $t \neq 0$, comparing the two sides yields
	\begin{equation}\label{eq:intermediate}
		g h_1 \alpha(h_2) = \alpha(h_1) S^{4}(h_2) g.
	\end{equation}
	Consequently,
	\[ S^4(h) = \alpha^{-1}(h_1) \alpha(h_2) S^{4}(h_3) g g^{-1} = \alpha^{-1}(h_1) g h_2 \alpha(h_3) g^{-1} = g (\alpha \rhu h \lhu \alpha^{-1}) g^{-1}, \]
	where the second equality follows from \eqref{eq:intermediate}.
	Thus, we obtain \eqref{eq:main-formula}.
\end{proof}

More generally, the same argument applies whenever $\Ext^i_H(\Bbbk, H)$ is one-dimensional and either the coradical of $H$ is coreflexive or the trivial left $H$-module is of type $\mathrm{FP}_{i}$.

In \Cref{thm:main}, the group-like element $g$ is uniquely determined by the one-dimensional left $H$-comodule structure of $\lint$, and the character $\alpha$ is uniquely determined by the one-dimensional right $H$-module structure of $\lint$.

We conclude this subsection by showing that, in the finite-dimensional case, the group-like element $g$ above recovers the classical distinguished group-like element.

\begin{rem}
	Let $H$ be a finite-dimensional Hopf algebra, and let $t$ be a nonzero left integral of $H$.
	Let $g$ denote the distinguished group-like element of $H$, characterized by the identity $\langle f, g \rangle \lambda = f * \lambda$ for any right integral $\lambda \in H^*$ and any $f \in H^*$.
	By \cite[Theorem 10.5.4]{Rad12}, $t_2 \otimes t_1 = t_1 \otimes S^2(t_2) g$.
	Applying $\big( (\id_H \otimes S^3) \circ \Delta \big) \otimes \id_H$ to both sides yields
	\[ t_2 \otimes S^3(t_3) \otimes t_1 = t_1 \otimes S^{3}(t_2) \otimes S^2(t_3) g. \]
	Multiplying the last two tensor factors yields
	\[ t_2 \otimes S^3(t_3) t_1 = t_1 \otimes S^{3}(t_2) S^2(t_3) g = t \otimes g. \]
	This shows that $g$ coincides with the group-like element of $H$ defined via the left coaction on $\lint$ in \Cref{thm:main}.
	Thus, our $S^4$-formula \eqref{eq:main-formula} agrees with Radford's formula when the Hopf algebra $H$ is finite-dimensional.
\end{rem}

\subsection{\texorpdfstring{$S^4$-formula for AS-co-Gorenstein Hopf algebras}{S4 formula for AS-co-Gorenstein Hopf algebras}}

In this subsection, we derive a formula for the fourth power of the antipode of an AS-co-Gorenstein Hopf algebra in terms of 
its homological cointegral.
The strategy is to equip the trivial comodule $\Bbbk^{\eta}$ with a coresolution by cofree $H$-comodules in a suitable category of Yetter--Drinfeld modules, so that the cohomology $\Ext_{\cM^H}^*(H, \Bbbk^{\eta})$ inherits a canonical twisted Yetter--Drinfeld structure.
We first construct a functor $- \boxplus H$ from right $H$-modules to Yetter--Drinfeld modules and show that it is right adjoint to the forgetful functor (\Cref{lem:right-adjoint}), thereby obtaining the desired $\cG$-coresolutions (\Cref{prop:coresolution}).
We then equip the relevant Ext groups with 
compatible $H$- and $H^*$-module structures that are independent of the choice of $\cG$-coresolution (\Cref{prop:Ext-structure-co}).
Under the rationality hypotheses, these structures yield twisted Yetter--Drinfeld modules.
Combining these ingredients, we prove the main result of this subsection (\Cref{thm:main-co}).

Let $M$ be a right $H$-module.
Equip $M \otimes H$ with the right $H$-comodule structure induced by the comultiplication on the second factor.
The linear map
\[ M \otimes H \otimes H \To M \otimes H, \qquad m \otimes x \otimes h \mapsto (m \otimes x) \lhd h := m \lhd h_2 \otimes \tau(S(h_1)) x \sigma(h_3) \]
endows $M \otimes H$ with a right $H$-module structure.
The module and comodule structures on $M \otimes H$ satisfy the compatibility condition \eqref{eq:yd-rr}.
Thus, $M \otimes H$ is a $(\sigma, \tau)$-twisted right-right Yetter--Drinfeld module, which we denote by $M \boxplus H$.
This construction yields a functor
\[- \boxplus H \colon \cM_H \To {}\!\YDrr(\sigma, \tau), \quad M \mapsto M \boxplus H. \]

This functor admits a left adjoint given by forgetting the comodule structure.

\begin{lem}\cite[Corollary 2.8]{CMZ97}\label{lem:right-adjoint}
	The functor $- \boxplus H \colon \cM_H \to \YDrr(\sigma, \tau)$ is right adjoint to the forgetful functor $\YDrr(\sigma, \tau) \to \cM_H$; that is, there is a natural isomorphism
	\[ \Psi_{V, M}: \Hom_{\cM_H}(V, M) \To \Hom_{\YDrr(\sigma, \tau)}(V, M \boxplus H), \qquad f \mapsto \big( v \mapsto f(v_0) \otimes v_1 \big). \]
\end{lem}

Let $\cG$ denote the class of objects in $\YDrr(\sigma,\tau)$ of the form $M \boxplus H$, where $M$ is a right $H$-module.

\begin{prop}\label{prop:coresolution}
	Let $V \in {\YDrr}(\sigma, \tau)$.
	\begin{enumerate}[label=\textup{(\arabic*)},leftmargin=2.2em]
		\item The morphism $\psi_V := \Psi_{V, V}(\id_V)$ is a \emph{$\cG$-preenvelope} of $V$; that is, for every $G \in \cG$, the induced map
		\[ (\psi_V)^* \colon \Hom_{{\YDrr}(\sigma, \tau)}(V \boxplus H, G) \To \Hom_{{\YDrr}(\sigma, \tau)}(V, G) \]
		is surjective.
		\item There exists a \emph{$\cG$-coresolution} of $V$; that is, a (not necessarily exact) complex $0 \to V \to G^0 \to G^1 \to \cdots$ in ${\YDrr}(\sigma,\tau)$ such that each $G^i \in \cG$ and the complex becomes exact upon applying $\Hom_{{\YDrr}(\sigma,\tau)}(-, G)$ for any $G \in \cG$.
		\item In particular, the $\cG$-coresolution of $V$ can be chosen to be an exact sequence
		\[ 0 \To V \stackrel{d^0}{\To} M^0 \boxplus H \stackrel{d^1}{\To} M^1 \boxplus H \stackrel{d^2}{\To} \cdots \stackrel{d^{i}}{\To} M^i \boxplus H \stackrel{d^{i+1}}{\To} \cdots \]
		in $\YDrr(\sigma, \tau)$, where, for any $i \geqslant 0$, $M^i$ is a right $H$-module and $M^i \cong \mathrm{Im}(d^i)$ as $H$-modules.
	\end{enumerate}
\end{prop}
\begin{proof}
	(1)
	Since $G \in \cG$, there exists an $H$-module $M$ such that $G = M \boxplus H$.
	For any $g \in \Hom_{{\YDrr}(\sigma,\tau)}(V, G)$, define $\widetilde{g} \in \Hom_{{\YDrr}(\sigma,\tau)}(V \boxplus H, G)$ by $\widetilde{g}(v \boxplus h) = \bigl( (\id_M \otimes \vep)g(v) \bigr) \otimes h$.
	For any $v \in V$,
	\[ (\widetilde{g} \psi_V)(v) = \widetilde{g}(v_0 \otimes v_1) = (\id_M \otimes \vep) g(v_0) \otimes v_1 = g(v), \]
	where the last equality follows from the $H$-colinearity of $g$.
	Therefore, $\widetilde{g} \psi_V = g$ and $(\psi_V)^*$ is surjective.
	
	(2) This is \cite[Proposition~8.1.3]{EJ11}.
	
	(3)
	Set $M^0 = V$ and $d^0 := \psi_{M^0}$. Since $d^0$ is injective, $M^0$ can be identified with $\mathrm{Im}(d^0)$.
	Set $M^1 := \Coker(d^0)$.
	Applying the same construction to $M^1$ and proceeding inductively, we obtain the desired exact sequence.
\end{proof}

By this adjunction, every twisted Yetter--Drinfeld module admits a coresolution by modules of the form $M \boxplus H$, which are cofree as $H$-comodules and hence injective.
Since each $M^i \boxplus H$ is injective as a right $H$-comodule, such a coresolution can be used to compute $\Ext$-groups of $H$-comodules.
To equip these $\Ext$-groups with a twisted Yetter--Drinfeld module structure, we first study the Hom spaces $\Hom_{\cM^H}(X, V)$ in a general setting.

\begin{lem}\label{lem:Hom-structure-co}
	Let $V \in \YDrr(\sigma, \tau)$ and $X \in {_H^H\!\cM^H_H}$.
	\begin{enumerate}[label=\textup{(\arabic*)},leftmargin=2.2em]
		\item The Hom space $\Hom_{\cM^H}(X, V)$ carries a natural left $H$-module structure given by
		\begin{equation}\label{eq:co-Hom-H-action}
			(h \rhd f)(x) = f\Bigl( \sigma\bigl( S^3(h_3) \bigr) \rhd x \lhd \tau(h_1) \Bigr) \lhd S(h_2).
		\end{equation}
		\item $\Hom_{\cM^H}(X, V)$ also carries a natural left $H^*$-module structure via
		\begin{equation}\label{eq:co-Hom-H^*-action}
			(h^* \rhd f)(x) = \langle h^*, x_{-1} \rangle f(x_0),
		\end{equation}
		and the resulting $H$-module and $H^*$-module structures on $\Hom_{\cM^H}(X, V)$ satisfy \eqref{eq:yd-dual-lr} for the pair $(\tau, \sigma S^4)$ of Hopf algebra endomorphisms.
		\item If the left $H^*$-module $\Hom_{\cM^H}(X, V)$ defined by \eqref{eq:co-Hom-H^*-action} is rational, then $\Hom_{\cM^H}(X, V) \in \YDlr(\tau, \sigma S^4)$.
	\end{enumerate}
\end{lem}
\begin{proof}
	(1)
    Obviously, the formula \eqref{eq:co-Hom-H-action} defines a left $H$-module structure on $\Hom(X, V)$.
	Thus, it remains to show that $(h \rhd f)$ is $H$-colinear for any $f \in \Hom_{\cM^H}(X, V)$ and $h \in H$.
	Using the compatibility condition and the antipode axioms, one obtains:
	\begin{align*}
		&  \rho\bigl( (h \rhd f)(x) \bigr) \\
		& = \rho\Bigl( f\bigl( \sigma(S^3(h_3)) \rhd x \lhd \tau(h_1) \bigr) \lhd S(h_2) \Bigr) \\
		& \xlongequal{\eqref{eq:yd-rr}} f\Bigl( \sigma\bigl( S^3(h_5) \bigr) \rhd x \lhd \tau(h_1) \Bigr)_0 \lhd S(h_3) \otimes \sigma\bigl( S^2(h_4) \bigr) f\Bigl( \sigma\bigl( S^3(h_5) \bigr) \rhd x \lhd \tau(h_1) \Bigr)_1 \tau\bigl( S(h_2) \bigr) \\
		& = f\Bigl( \sigma\bigl( S^3(h_7) \bigr) \rhd x_0 \lhd \tau(h_1) \Bigr) \lhd S(h_4) \otimes \sigma\bigl( S^2(h_5) S^3(h_6) \bigr) x_1 \tau\bigl( h_2 S(h_3) \bigr) \\
		& = f\Bigl( \sigma\bigl( S^3(h_3) \bigr) \rhd x_0 \lhd \tau(h_1) \Bigr) \lhd S(h_2) \otimes x_1 \\
		& = (h \rhd f)(x_0) \otimes x_1.
	\end{align*}
	
	(2)
	Formula \eqref{eq:co-Hom-H^*-action} endows $\Hom_{\cM^H}(X, V)$ with a well-defined left $H^*$-module structure induced by the left $H$-comodule structure on $X$.
	Since 
	\begin{align*}
		&  \Bigl( h_1 \rhd \bigl( (h^* \lhu \sigma(S^4(h_2)) ) \rhd f \bigr) \Bigr) (x) \\
		& = \Bigl( \bigl( h^* \lhu \sigma(S^4(h_4)) \bigr) \rhd f \Bigr) \Bigl( \sigma\bigl( S^3(h_3) \bigr) \rhd x \lhd \tau(h_1) \Bigr) \lhd S(h_2) \\
		& = f\Bigl( \bigl( \sigma(S^3(h_3)) \rhd x \lhd \tau(h_1) \bigr) \lhd \bigl(h^* \lhu \sigma(S^4(h_4)) \bigr) \Bigr) \lhd S(h_2) \\
		& = \bigl\langle h^* \lhu \sigma(S^4(h_6)), \sigma(S^3(h_5)) x_{-1} \tau(h_1) \bigr\rangle f \Bigl( \sigma\bigl( S^3(h_4) \bigr) \rhd x_0 \lhd \tau(h_2) \Bigr) \lhd S(h_3) \\
		& = \bigl\langle h^*, \sigma\bigl( S^4(h_6) S^3(h_5) \bigr) x_{-1} \tau(h_1) \bigr\rangle f \Bigl( \sigma\bigl( S^3(h_4) \bigr) \rhd x_0 \lhd \tau(h_2) \Bigr) \lhd S(h_3) \\
		& = \bigl\langle h^*, x_{-1} \tau(h_1) \bigr\rangle f \Bigl( \sigma\bigl( S^3(h_4) \bigr) \rhd x_0 \lhd \tau(h_2) \Bigr) \lhd S(h_3) \\
		& = \langle \tau(h_1) \rhu h^*, x_{-1} \rangle (h_2 \rhd f)(x) \\
		& = \Bigl( \bigl( \tau(h_1) \rhu h^* \bigr) \rhd (h_2 \rhd f) \Bigr)(x),
	\end{align*}
	it follows that \eqref{eq:yd-dual-lr} holds for the pair $(\tau, \sigma S^4)$ of Hopf algebra endomorphisms.
	
	(3)
	This follows from (2) and the equivalence of \eqref{eq:yd-lr} and \eqref{eq:yd-dual-lr} when the left $H^*$-module structure is induced by a right $H$-comodule structure.
\end{proof}

We use the following comparison lemma, whose proof is routine and is therefore omitted (cf. \cite[Section~8.1]{EJ11}).

\begin{lem}\label{lem:comparison-co}
	Let $U, V \in \YDrr(\sigma, \tau)$, let $0 \to U \to G^0 \to G^1 \to \cdots$ be a $\cG$-coresolution of $U$, and let $0 \to V \to Q^0 \to Q^1 \to \cdots$ be a complex in $\YDrr(\sigma, \tau)$ with $Q^i \in \cG$ for all $i \geqslant 0$.
	Then every morphism $f \colon U \to V$ in $\YDrr(\sigma, \tau)$ lifts to a cochain map $G^{\bullet} \to Q^{\bullet}$, unique up to homotopy.
\end{lem}

By \Cref{lem:comparison-co}, every morphism in $\YDrr(\sigma, \tau)$ lifts to a cochain map of $\cG$-coresolutions, and the lift is unique up to homotopy.
The following proposition uses this lift.

\begin{prop}\label{prop:Ext-structure-co}
	Let $V \in \YDrr(\sigma, \tau)$ and let $X \in {_H^H\!\cM^H_H}$ be a Hopf bimodule.
	For every $i \geqslant 0$, the following statements hold.
	\begin{enumerate}[label=\textup{(\arabic*)},leftmargin=2.2em]
		\item The space $\Ext_{\cM^H}^i(X, V)$ carries a natural left $H^*$-module structure induced by \eqref{eq:co-Hom-H^*-action} and a natural left $H$-module structure induced by \eqref{eq:co-Hom-H-action}.
		Moreover, both structures are independent of the choice of $\cG$-coresolution of $V$ constructed in \Cref{prop:coresolution}.
		\item The space $\Ext_{\cM^H}^i(X, V)$ belongs to $\widehat{\YDlr}(\tau, \sigma S^4)$.
		If this space is also rational as a left $H^*$-module, then it belongs to $\YDlr(\tau, \sigma S^4)$.
	\end{enumerate}
\end{prop}
\begin{proof}
	(1)
	Let $G^{\bullet}$ be the $\cG$-coresolution of $V$ given in \Cref{prop:coresolution}(3).
	This coresolution is exact, and each of its terms is injective as a right $H$-comodule.
	Hence $\mathrm{H}^i\bigl( \Hom_{\cM^H}(X, G^{\bullet}) \bigr)$ is isomorphic, as a $\Bbbk$-vector space, to the usual Ext group $\Ext_{\cM^H}^i(X, V)$ computed in the category of right $H$-comodules.
	By \Cref{lem:Hom-structure-co}, the complex $\Hom_{\cM^H}(X, G^{\bullet})$ is both an $H$-complex and an $H^*$-complex.
	Hence, its cohomology groups inherit left $H$-module and left $H^*$-module structures.

	Suppose that ${G'}^{\bullet}$ is another $\cG$-coresolution of $V$.
	By \Cref{lem:comparison-co}, there exists a cochain map $f^{\bullet} \colon G^{\bullet} \to {G'}^{\bullet}$ that lifts $\id_V$ and is a homotopy equivalence.
	Thus
	\[ \Hom_{\cM^H}(X, f^{\bullet}) \colon \Hom_{\cM^H}(X, G^{\bullet}) \To \Hom_{\cM^H}(X, {G'}^{\bullet}) \]
	is also a homotopy equivalence.
	Since each $\Hom_{\cM^H}(X, f^i)$ is both $H$-linear and $H^*$-linear, the two induced module structures on the cohomology groups agree.

	(2)
	By \Cref{lem:Hom-structure-co}(2), each term of $\Hom_{\cM^H}(X, G^{\bullet})$ satisfies \eqref{eq:yd-dual-lr} for the pair $(\tau, \sigma S^4)$. The differentials are morphisms of both module structures, so the identity passes to the cohomology, and $\Ext_{\cM^H}^i(X, V)$ belongs to $\widehat{\YDlr}(\tau, \sigma S^4)$. If this left $H^*$-module is rational, then $\Ext_{\cM^H}^i(X, V)$ belongs to $\YDlr(\tau, \sigma S^4)$, as in \Cref{lem:Hom-structure-co}(3).
\end{proof}

We now establish an analogue of \Cref{thm:main} for AS-co-Gorenstein Hopf algebras.

\begin{thm}\label{thm:main-co}
	Let $H$ be a right AS-co-Gorenstein Hopf algebra of dimension $d$ with antipode $S$.
	Suppose that the coradical of $H$ is coreflexive.
	\begin{enumerate}[label=\textup{(\arabic*)},leftmargin=2.2em]
		\item The right homological cointegral space $\rcoint$ is an $(\id_H, S^4)$-twisted left-right Yetter--Drinfeld module.
		Thus, there exist a group-like element $\widetilde{g}$ and a character $\widetilde{\alpha}$ such that $\rcoint \cong \Bbbk{\widetilde{g}}$ as a right $H$-comodule and $\rcoint \cong {}_{\widetilde{\alpha}}\Bbbk$ as a left $H$-module.
		\item For all $h\in H$,
		\begin{equation}\label{eq:main-co-formula}
			S^4(h) = \widetilde{g} (\widetilde{\alpha} \rhu h \lhu \widetilde{\alpha}^{-1}) \widetilde{g}^{-1}.
		\end{equation}
	\end{enumerate}
\end{thm}
\begin{proof}
	(1)
	Note that $\Bbbk$, equipped with the trivial right $H$-module and right $H$-comodule structures, is a right-right Yetter--Drinfeld module.
	Applying \Cref{prop:Ext-structure-co} with $\sigma = \tau = \id_H$, we obtain
	\[ \Ext^i_{\cM^H}(H, \Bbbk) \in \widehat{\YDlr}(\id_H, S^4) \quad \text{for all } i \geqslant 0. \]
    Since $H$ is right AS-co-Gorenstein of dimension $d$, the space $\rcoint = \Ext^d_{\cM^H}(H, \Bbbk)$ is one-dimensional and hence simple.
    Moreover, $\rcoint$ is rational by \Cref{lem:coreflexive-rational}, and therefore $\rcoint \in \YDlr(\id_H, S^4)$.
	Consequently, there exist a group-like element $\widetilde{g} \in G(H)$ and a character $\widetilde{\alpha} \colon H \to \Bbbk$ such that
    $\rcoint \cong \Bbbk{\widetilde{g}}$ as a right $H$-comodule and $\rcoint \cong {}_{\widetilde{\alpha}}\Bbbk$ as a left $H$-module.

    (2)
    Fix a nonzero element $t \in \rcoint$.
    By the compatibility condition \eqref{eq:yd-lr} with $(\sigma, \tau) = (\id_H, S^4)$, we obtain
    \[ \widetilde{\alpha}(h_1)\, t \otimes S^4(h_2)\widetilde{g}
    = (h_2 \rhd t)_0 \otimes (h_2 \rhd t)_1 h_1
    = \widetilde{\alpha}(h_2)\, t \otimes \widetilde{g}\, h_1. \]
    Since $t \neq 0$, comparing the two sides yields
    \[ \widetilde{g}\, h_1\, \widetilde{\alpha}(h_2) = \widetilde{\alpha}(h_1)\, S^4(h_2)\, \widetilde{g}. \]
    This is \eqref{eq:intermediate} with $(g,\alpha)$ replaced by $(\widetilde{g},\widetilde{\alpha})$. The remainder of the proof of \Cref{thm:main}(2) therefore yields \eqref{eq:main-co-formula}.
\end{proof}

If $H$ is a right AS-co-Gorenstein Hopf algebra of dimension $0$, that is, $H$ is co-Frobenius, then $\rcoint$ is a one-dimensional $H^*$-submodule of $H^*$.
Since every finite-dimensional $H^*$-submodule of $H^*$ is rational \cite[Proposition~3.1.9(b)]{Rad12}, $\rcoint$ is rational and hence an $(\id_H, S^{4})$-twisted left-right Yetter--Drinfeld module.
In this case, $\widetilde{g}$ and $\widetilde{\alpha}$ coincide with the distinguished group-like elements in \cite{BBT07}; see (1.8) and Lemma~2.5 therein.
As in the proof of \Cref{thm:main}, 
the $S^{4}$ formula for co-Frobenius Hopf algebras \cite[Theorem~2.8]{BBT07} follows directly from the twisted Yetter--Drinfeld compatibility condition.

\section{\texorpdfstring{Applications of the $S^4$-formula}{Applications of the S4 formula}}\label{sec:app}

In this section, we present several applications of the $S^4$-formulas established in \Cref{thm:main} and \Cref{thm:main-co}.
We first show, in \Cref{prop:trivial-ASG} for homological integrals and in \Cref{prop:trivial-AScoG} for homological cointegrals, that the module structure is trivial when $H$ is commutative and the comodule structure is trivial when $H$ is cocommutative.
Group algebras provide a special case (\Cref{cor:kG-cointegral}).
We then compare $(\lint)^*$ with $\rcoint$ for the fundamental group of the Klein bottle and for the quantized Borel subalgebra ${\cal U}^{\geqslant 0}_q(\mathfrak{sl}_2)$.
In the final two subsections, we apply the $S^4$-formulas to two structural questions: the bijectivity of the antipode and the finiteness of the order of $S^2$.

\subsection{\texorpdfstring{Triviality of homological (co-)integrals and comparison of $(\lint)^*$ with $\rcoint$}{Triviality of homological (co-)integrals and comparison of integrals with cointegrals}}

The next two propositions separate the two structures: if $H$ is commutative, then the module structures on its homological integrals and cointegrals are trivial, while if $H$ is cocommutative, then the comodule structures are trivial.

\begin{prop}\label{prop:trivial-ASG}
	Let $H$ be an AS-Gorenstein Hopf algebra.
	\begin{enumerate}[label=\textup{(\arabic*)},leftmargin=2.2em]
		\item If $H$ is commutative, then $\lint \cong \Bbbk_{\vep}$ as a right $H$-module.
		\item If $H$ is cocommutative, then $\lint$ is rational as a right $H^*$-module, and $\lint \cong {}^{\eta}\Bbbk$ as a left $H$-comodule.
	\end{enumerate}
\end{prop}
\begin{proof}
	(1)
	For $z\in H$, let $r_z \colon H \to H$, $x \mapsto xz$, denote right multiplication by $z$. By definition, the $H$-module structure on $\lint = \Ext^d_H(\Bbbk, H)$ is induced by $r_z$.
	If $H$ is commutative and $z \in H$, then $r_z$ coincides with left multiplication by $z$, which is $H$-linear.
	Multiplication by $z$ on any projective resolution of $\Bbbk$ is a chain map lifting $\vep(z) \id_{\Bbbk}$. Hence, $r_z$ induces $\vep(z) \id$ on $\Ext^d_H(\Bbbk, H)$.
	Consequently, the $H$-module structure on $\lint$ factors through the counit, i.e.\ $\lint \cong \Bbbk_{\vep}$ as $H$-modules.
	
	(2)
	For any left $H$-comodule $N$ whose coaction is induced by the unit $1_H$, that is, $n_{-1} \otimes n_0 = 1_H \otimes n$, the coaction on $H \boxtimes N$ satisfies
	\[ \rho(h \boxtimes n) = h_1 n_{-1} S(h_3) \otimes h_2 \boxtimes n_0 = h_1 1_H S(h_3) \otimes h_2 \boxtimes n = h_1 S(h_2) \otimes h_3 \boxtimes n = 1_H \otimes h \boxtimes n, \]
	where the third equality follows from the cocommutativity of $H$.
	In particular, the coaction on $H \boxtimes N$ is trivial.
	Now consider the exact sequence
	\[ \cdots \stackrel{d_{i+1}}{\To} H \boxtimes N_i \stackrel{d_i}{\To} \cdots \stackrel{d_2}{\To} H \boxtimes N_1 \stackrel{d_1}{\To} H \boxtimes N_0 \stackrel{d_0}{\To} \Bbbk \To 0 \]
	from \Cref{prop:resolution}, where each $d_i|_{1 \boxtimes N_i}$ is injective.
	Since $d_0|_{1 \boxtimes N_0}$ embeds $N_0$ as a subcomodule of $\Bbbk$, $N_0$ is isomorphic to $\Bbbk$ as a left $H$-comodule.
	Thus, the coaction on $H \boxtimes N_0$ is trivial.
	Since $d_1|_{1 \boxtimes N_1}$ is an injective colinear map into $H \boxtimes N_0$, the coaction on $N_1$ is also trivial.
	Repeating this argument, we conclude by induction that the coaction on $H \boxtimes N_i$ is trivial for every $i \geqslant 0$.
	For any $h \in H$, $h^* \in H^*$, $n \in N_i$, and $f \in \Hom_H(H \boxtimes N_i, H)$, we have
	\[ (f \lhd h^*)(h \boxtimes n) \xlongequal{\eqref{eq:Hom-H^*-action}} f(h \otimes n)_2 \bigl\langle h^*, S^3\bigl( f(h \boxtimes n)_{3} \bigr) S(1_H) f(h \boxtimes n)_{1} \bigr\rangle = \langle h^*, 1 \rangle f(h \boxtimes n). \]
    Thus, the $H^*$-action on $\Hom_H(H \boxtimes N_i, H)$ factors through evaluation at the unit $1_H$.
	Since this holds for every $i \geqslant 0$, the $H^*$-action on the entire complex $\Hom_H(H \boxtimes N_\bullet, H)$, and hence on its cohomology, factors through evaluation at $1_H$.
	Therefore, $\lint \cong {}^{\eta}\Bbbk$ as left $H$-comodules.
\end{proof}

The following proposition is the co-Gorenstein analogue of \Cref{prop:trivial-ASG}; its proof uses the coresolution in 
\Cref{prop:coresolution} in place of the resolution in \Cref{prop:resolution}.

\begin{prop}\label{prop:trivial-AScoG}
	Let $H$ be an AS-co-Gorenstein Hopf algebra.
	\begin{enumerate}[label=\textup{(\arabic*)},leftmargin=2.2em]
		\item If $H$ is commutative, then $\rcoint \cong {}_{\vep}\Bbbk$ as a left $H$-module.
		\item If $H$ is cocommutative, then $\rcoint$ is rational as a left $H^*$-module, and $\rcoint \cong \Bbbk^{\eta}$ as a right $H$-comodule.
	\end{enumerate}
\end{prop}
\begin{proof}
	(1)
	Applying \Cref{prop:coresolution} to $V = \Bbbk_{\vep}^{\eta} \in \YDrr$, we obtain an exact sequence
	\[ 0 \To \Bbbk \stackrel{d^0}{\To} M^0 \boxplus H \stackrel{d^1}{\To} M^1 \boxplus H \stackrel{d^2}{\To} \cdots \stackrel{d^{i}}{\To} M^i \boxplus H \stackrel{d^{i+1}}{\To} \cdots \]
	in $\YDrr$, where $M^i \cong \mathrm{Im}(d^i)$ as $H$-modules for every $i \geqslant 0$.
	For any right $H$-module $M$ with trivial $H$-action, 
	the right $H$-action on $M \boxplus H$ satisfies
	\[ (m \boxplus k) \lhd h = m \lhd h_2 \boxplus S(h_1) k h_3 = m \boxplus S(h_1) k h_2 = m \boxplus k S(h_1) h_2 = \vep(h) (m \boxplus k), \]
	where the third equality follows from the commutativity of $H$.
	As shown in the proof of \Cref{prop:coresolution}, $M^0$ is isomorphic to $\Bbbk$, so $M^0 \boxplus H$ is a trivial $H$-module. Consequently, its quotient $\mathrm{Im}(d^1)$ is also a trivial $H$-module.
	Moreover, since each $M^{i+1} \cong \mathrm{Im}(d^{i+1})$ as $H$-modules, it follows by induction that $\mathrm{Im}(d^{i+1})$, $M^{i+1}$, and $M^{i+1} \boxplus H$ are trivial for every $i \geqslant 0$.
	Thus, the $H$-action on each $M^i \boxplus H$ is trivial. 
	Moreover, for any $f \in \Hom_{\cM^H}(H, M^i \boxplus H)$, the action \eqref{eq:co-Hom-H-action} yields
	\[ (h \rhd f)(k) = f\bigl( S^3(h_3) k h_1 \bigr) \lhd S(h_2) = f\bigl( S^3(h_2) k h_1 \bigr) = f\bigl( k h_1 S^3(h_2) \bigr) = \vep(h) f(k), \]
	so the $H$-action on each $\Hom_{\cM^H}(H, M^i \boxplus H)$ is induced by the counit $\vep$.
	Therefore, the $H$-action on the entire complex $\Hom_{\cM^H}(H, M^{\bullet} \boxplus H)$, and hence on its cohomology $\Ext^i_{\cM^H}(H, \Bbbk)$, factors through the counit.
	We conclude that $\rcoint \cong {}_{\vep}\Bbbk$ as $H$-modules.
	
	(2)
	Applying \Cref{prop:coresolution} to $V = \Bbbk_{\vep}^{\eta} \in \YDrr$, we obtain an exact sequence
	\[ 0 \To \Bbbk_{\vep}^{\eta} \stackrel{d^0}{\To} G^0 \stackrel{d^1}{\To} G^1 \stackrel{d^2}{\To} \cdots \stackrel{d^{i}}{\To} G^i \stackrel{d^{i+1}}{\To} \cdots \]
	in $\YDrr$, where each $G^i$ is an injective $H$-comodule.
	For any $\lambda \in \Hom_{\cM^H}(H, G^i)$, the $H$-colinearity of $\lambda$ and the cocommutativity of $H$ imply that
	\[ \lambda(h_2) \otimes h_1 = \lambda(h_1) \otimes h_2 = \lambda(h)_0 \otimes \lambda(h)_1. \]
	Therefore, for $h^* \in H^*$, the action \eqref{eq:co-Hom-H^*-action} yields
	\[ (h^* \rhd \lambda)(h) = \langle h^*, h_1 \rangle \lambda(h_2) =  \langle h^*, \lambda(h)_1 \rangle \lambda(h)_0 = h^* \rhd \lambda(h). \]
	Thus, the $H^*$-action \eqref{eq:co-Hom-H^*-action} on each $\Hom_{\cM^H}(H, G^i)$ agrees with the action induced by the $H^*$-action on $G^i$.
	Since $H$ is cocommutative, the action map $h^* \rhd -$ is right $H$-colinear on each $G^i$ and thus defines an $H$-colinear chain endomorphism $\xi$ of $G^{\bullet}$.
	Since $h^* \rhd -$ coincides with $\langle h^*, 1_H \rangle \id$ on $\Bbbk_{\vep}^{\eta}$, $\xi \colon G^{\bullet} \to G^{\bullet}$ is homotopic to $\langle h^*, 1_H \rangle \id_{G^{\bullet}}$.
	Applying $\Hom_{\cM^H}(H, -)$, we conclude that the induced chain endomorphism $\Hom_{\cM^H}(H, \xi)$ of $\Hom_{\cM^H}(H, G^{\bullet})$ is likewise homotopic to $\langle h^*, 1_H \rangle \id$.
	Since the chain endomorphism defined by the action of $h^*$ on $\Hom_{\cM^H}(H, G^{\bullet})$ coincides with $\Hom_{\cM^H}(H, \xi)$, it is homotopic to $\langle h^*, 1_H \rangle \id$.
	Therefore, the induced action of $h^*$ on $\Ext_{\cM^H}^i(H, \Bbbk^{\eta})$ is given by the scalar $\langle h^*, 1_H \rangle$; equivalently, the $H^*$-action on $\Ext_{\cM^H}^i(H, \Bbbk^{\eta})$ is induced by a trivial $H$-coaction.
	By the AS-co-Gorenstein hypothesis, $\rcoint$ is one-dimensional.
	Thus, $\rcoint \cong \Bbbk^{\eta}$ as $H$-comodules.
\end{proof}

\Cref{prop:trivial-AScoG}(2) gives the comodule statement in the following corollary. 

\begin{cor}\label{cor:kG-cointegral}
	Let $G$ be a group.
	Then $\textstyle \int_{r}^{\Bbbk G}$ is rational as a left $(\Bbbk G)^*$-module, and $\textstyle \int_{r}^{\Bbbk G} \cong {}_{\vep}\Bbbk^{\eta}$ as Yetter--Drinfeld modules.
\end{cor}
\begin{proof}
	As a coalgebra, the group algebra $\Bbbk G$ is cocommutative and decomposes as a direct sum of the one-dimensional subcoalgebras $\Bbbk g$ ($g\in G$).
	Let $\lambda \in \Hom(\Bbbk G,\Bbbk)$ be the linear map defined by
	\[ \lambda(g) = \begin{cases}
		1, & g = 1_G, \\
		0, & g \neq 1_G.
	\end{cases} \]
	Since $\Delta(g) = g \otimes g$, one has $(\lambda * f)(g) = \lambda(g)\, f(g) = \langle f, 1_G \rangle \lambda(g)$ for every $f \in (\Bbbk G)^*$ and every $g \in G$. Thus $\lambda * f = \langle f, 1_G \rangle \lambda$, so $\Bbbk G$ is co-Frobenius and hence AS-co-Gorenstein of dimension $0$. In particular, $\textstyle\int_{r}^{\Bbbk G} = \Hom_{\cM^{\Bbbk G}}(\Bbbk G, \Bbbk)$.
	The map $\lambda$ is $\Bbbk G$-colinear.
	For any $g, h \in G$, the module action in \Cref{lem:Hom-structure-co} yields
	\[ (g \rhd \lambda)(h) = g \rhd \lambda(g^{-1} h g) = \lambda(g^{-1}hg) = \begin{cases}
		1, & h = 1_G, \\
		0, & h \neq 1_G.
	\end{cases} \]
	Hence, $g \rhd \lambda = \lambda$, and therefore $\textstyle \int_{r}^{\Bbbk G} \cong {}_{\vep} \Bbbk$ as $\Bbbk G$-modules.
    Together with \Cref{prop:trivial-AScoG}(2), this shows that $\textstyle \int_{r}^{\Bbbk G} \cong {}_{\vep}\Bbbk^{\eta}$ as Yetter--Drinfeld modules.
\end{proof}

When $H$ is finite-dimensional, one can show that $(\lint)^* \cong \rcoint$ as $(\id_H, S^4)$-twisted left-right Yetter--Drinfeld modules, where the twisted Yetter--Drinfeld module structure on $(\lint)^*$ is given by \Cref{lem:dual-YD}.
The following example shows that this duality does not hold in general, even for noetherian AS-regular group algebras.

\begin{ex}[The Klein bottle]\label{ex:klein}
	Assume that $\cha(\Bbbk) \neq 2$. Let $K$ be the Klein bottle, and set $G = \pi_1(K)$. Then
	\[ G \cong \langle a, b \mid aba^{-1} = b^{-1} \rangle \cong \mathbb{Z} \rtimes_{-1}\mathbb{Z}. \]
	Since $G$ is polycyclic, $\Bbbk G$ is a noetherian group algebra by \cite[Theorem 1.5.12]{MR01}.
	There is an exact sequence 
	\[ 0 \To \Bbbk G \boxtimes N_2 \stackrel{d_2}{\To} \Bbbk G \boxtimes N_1 \stackrel{d_1}{\To} \Bbbk G \boxtimes N_0 \stackrel{d_0}{\To} \Bbbk \To 0 \]
	in $^{\Bbbk G}_{\Bbbk G}\!\YD$, where $N_0 = \Bbbk e_0$, $N_1 = \Bbbk e_{11} \oplus \Bbbk e_{12}$, and $N_2 = \Bbbk e_2$; the $\Bbbk G$-coactions on the $N_i$ are induced by the unit $1_G$; and the boundary homomorphisms are determined by
	\[ d_0(1 \boxtimes e_0) = 1, \qquad d_1(1 \boxtimes e_{11}) = (a-1) \boxtimes e_{0}, \qquad d_1(1 \boxtimes e_{12}) = (b-1) \boxtimes e_0, \]
	\[ d_2(1 \boxtimes e_2) = (b-1) \boxtimes e_{11} + (ab^{-1}+1) \boxtimes e_{12}. \]
	Applying the functor $\Hom_{\Bbbk G}(-, \Bbbk G)$ to this projective resolution, we obtain the complex
	\[ 0 \To \Hom_{\Bbbk G}(\Bbbk G \boxtimes N_0, \Bbbk G) \stackrel{d_1^*}{\To} \Hom_{\Bbbk G}(\Bbbk G \boxtimes N_1, \Bbbk G) \stackrel{d_2^*}{\To} \Hom_{\Bbbk G}(\Bbbk G \boxtimes N_2, \Bbbk G) \To 0. \]
	A direct computation of the cohomology shows that
	\[ \Ext^i_{\Bbbk G}({}_{\vep}\Bbbk, \Bbbk G) \cong \begin{cases}
		0, & i \neq 2, \\
		\Bbbk G/ \bigl( (b-1)\Bbbk G + (ab^{-1}+1)\Bbbk G \bigr) \cong \Bbbk_{\alpha}, & i = 2,
	\end{cases} \]
	where $\alpha$ is the character of $\Bbbk G$ defined  by $\alpha(a) = -1$ and $\alpha(b) = 1$. 
	Hence $\Bbbk G$ is AS-regular of dimension $2$.
    By \Cref{prop:trivial-ASG}(2), $\textstyle \int^{\ell}_{\Bbbk G} \cong {}^{\eta}\Bbbk$ as comodules.
	Consequently, $(\textstyle \int^{\ell}_{\Bbbk G})^* \cong ({}^{\eta}\Bbbk_{\alpha})^* \cong {}_{\alpha}\Bbbk^{\eta} \ncong {}_{\vep}\Bbbk^{\eta}$, whereas $\textstyle \int_{r}^{\Bbbk G} \cong {}_{\vep}\Bbbk^{\eta}$ by \Cref{cor:kG-cointegral}.
	Thus, there exists a noetherian AS-regular group algebra $\Bbbk G$ such that $(\textstyle \int^{\ell}_{\Bbbk G})^* \ncong \textstyle \int_{r}^{\Bbbk G}$.
\end{ex}

In contrast to the preceding counterexample, the following example demonstrates that the duality $(\int_{H}^{\ell})^* \cong \int^{H}_r$ still holds for the quantized Borel subalgebra ${\cal U}^{\geqslant 0}_q(\mathfrak{sl}_2)$, which is AS-regular of dimension $2$ and AS-co-regular of dimension $1$.

\begin{ex}
	Assume that $q \in \Bbbk^{\times}$ is not a root of unity.
	The quantized Borel subalgebra $H:= {\cal U}^{\geqslant 0}_q(\mathfrak{sl}_2)$ is the $\Bbbk$-algebra generated by $E, K, K^{-1}$ subject to the relations
	\[ KK^{-1} = 1 = K^{-1}K, \qquad KEK^{-1} = q^2E. \]
	The coalgebra structure is given by
	\[ \Delta(K) = K \otimes K, \qquad \Delta(E) = E \otimes 1 + K \otimes E, \qquad \vep(K) = 1, \qquad \vep(E) = 0, \]
	and the antipode $S$ is given by $S(K) = K^{-1}$ and $S(E) = - K^{-1}E$.
	
	The algebra $H$ is noetherian AS-regular of dimension $2$.
	There is an exact sequence
	\[ 0 \To H \boxtimes N_2 \stackrel{d_2}{\To} H \boxtimes N_1 \stackrel{d_1}{\To} H \boxtimes N_0 \stackrel{d_0}{\To} \Bbbk \To 0 \]
	in ${}^{H}_{H}\YD$, where $N_0 = \Bbbk e_0$, $N_1 = \Bbbk e_{11} \oplus \Bbbk e_{12}$ and $N_2 = \Bbbk e_2$. The $H$-coactions on the $N_i$ are given by
	\[ \rho(e_0) = 1 \otimes e_0, \qquad \rho(e_{11}) = 1 \otimes e_{11}, \qquad \rho(e_{12}) = K \otimes e_{12} - E \otimes e_{11}, \qquad \rho(e_2) = K \otimes e_2, \]
	and the boundary homomorphisms are determined by
	\[ d_0(1 \boxtimes e_0) = 1, \qquad d_1(1 \boxtimes e_{11}) = (K-1) \boxtimes e_{0}, \qquad d_1(1 \boxtimes e_{12}) = E \boxtimes e_0, \]
	\[ d_2(1 \boxtimes e_2) = E \boxtimes e_{11} + (1 - q^{-2}K) \boxtimes e_{12}. \]
	Applying $\Hom_{H}(-, H)$ to this projective resolution yields the complex
	\[ 0 \To \Hom_{H}(H \boxtimes N_0, H) \stackrel{d_1^*}{\To} \Hom_{H}(H \boxtimes N_1, H) \stackrel{d_2^*}{\To} \Hom_{H}(H \boxtimes N_2, H) \To 0. \]
	A direct computation shows that
	\[ \Ext^i_{H}({}_{\vep}\Bbbk, H) = 0 \qquad \text{for all } i \neq 2. \]
	Moreover, $\textstyle \int^{\ell}_{H} = \Ext^2_{H}({}_{\vep}\Bbbk, H) \cong \Bbbk_{\alpha}$ as right $H$-modules, where $\alpha(K) = q^2$ and $\alpha(E) = 0$, and $\textstyle \int^{\ell}_{H} \cong \Bbbk K^{-1}$ as a left $H$-comodule.
	
	Next, we show that $H$ is AS-co-regular of dimension $1$ and compute the right homological cointegral $\textstyle \int^{H}_r$.
	Let $\Bbbk[x]$ be the subalgebra of $H$ generated by $x = K^{-1}E$.
	The subalgebra $\Bbbk[x]$ is a right coideal of $H$ and is stable under the right adjoint action $b \lhd h = S(h_1) b h_2$ of $H$.
	Since a Hopf algebra, equipped with the right adjoint action, is a right-right Yetter--Drinfeld module over itself, $\Bbbk[x]$ is a subobject of $H$ in $\YDrr$.
	Moreover, $H = \bigoplus_{n \in \Z} K^n\Bbbk[x]$ as right comodules, and $H$ is injective as a right comodule over itself; hence each direct summand $K^n\Bbbk[x]$ is injective.
	Let $\pi_n$ denote the canonical projection of $H$ onto $K^n\Bbbk[x]$.
	Let $p = q^2$.
	Set $k_p = 1 + p + \cdots + p^{k-1}$, $0!_p = 1$, and $k!_p = 1_p 2_p \cdots k_p$ for any $k \geqslant 1$.
	Define the Gaussian binomial coefficient by $\binom{n}{k}_p = \frac{n!_p}{k!_p (n-k)!_p}$ for $0 \leqslant k \leqslant n$.
	Equip $K^{-1}\Bbbk[x]$ with the right $H$-module structure given by
	\[ (K^{-1}x^n) \lhd K = q^{-2(n+1)} K^{-1}x^n, \qquad (K^{-1}x^n) \lhd E = (q^{-2(n+1)} - q^2) K^{-1}x^{n+1} \]
	for all $n \geqslant 0$.
	The right $H$-coaction on $K^{-1}\Bbbk[x]$ is the restriction of $\Delta$. Together with the right action above, it satisfies \eqref{eq:yd-rr}, so $K^{-1}\Bbbk[x]$ is an object of $\YDrr$.
	There is an exact sequence in $\YDrr$:
	\begin{equation}\label{eq:G-coresolution-borel}
		0 \To \Bbbk^{\eta}_{\vep} \stackrel{d^0}{\To} \Bbbk[x] \stackrel{d^1}{\To} K^{-1}\Bbbk[x] \To 0
	\end{equation}
	The boundary homomorphisms are given by
	\[ d^0(1) = 1, \qquad d^1(x^n) = \begin{cases}
		0, & n = 0, \\
		n_{q^{2}} K^{-1} x^{n-1}, & n \geqslant 1.
	\end{cases} \]
    Since \eqref{eq:G-coresolution-borel} is split exact as a sequence of right $H$-modules, applying $\Hom_{\cM_H}(-, M)$ to it yields an exact sequence for every right $H$-module $M$.
    By the natural isomorphism of \Cref{lem:right-adjoint}, $\Hom_{\YDrr}(-, M \boxplus H) \cong \Hom_{\cM_H}(-, M)$.
    Hence applying $\Hom_{\YDrr}(-, G)$ to \eqref{eq:G-coresolution-borel} likewise yields an exact sequence for every $G = M \boxplus H \in \cG$, and \eqref{eq:G-coresolution-borel} is a $\cG$-coresolution of $\Bbbk^{\eta}_{\vep}$ in $\YDrr$.
	Since it is also an injective resolution of $\Bbbk^{\eta}_{\vep}$ of length $1$, the injective dimension of the trivial $H$-comodule $\Bbbk$ is at most $1$.
	
	For any $f \in \Hom_{\cM^{H}}(H, \Bbbk[x])$,
	we have
	\[ f(K^mx^n) = \begin{cases}
		\binom{n}{n-m}_{q^{2}} \vep(f(K^mx^{m})) x^{n-m}, & 0 \leqslant m \leqslant n \\
		0, & \text{otherwise},
	\end{cases} \]
	and
	\[ (d^{1} \circ f)(K^mx^n) = \begin{cases}
		\binom{n}{n-m}_{q^{2}} (n-m)_{q^{2}} \vep(f(K^mx^{m})) K^{-1} x^{n-m-1}, & 0 \leqslant m < n \\
		0, & \text{otherwise}.
	\end{cases} \]
	If $d^{1} \circ f = 0$, then setting $n = m+1$ in the formula above yields $\vep(f(K^mx^{m})) = 0$ for all $m \geqslant 0$, and hence $f = 0$.
	Thus $\Ext_{\cM^{H}}^0(H, \Bbbk) = 0$.
	Moreover, $d^1 \circ f \neq \pi_{-1}$ for any $f \in \Hom_{\cM^{H}}(H, \Bbbk[x])$.
	On the other hand, for any $g \in \Hom_{\cM^{H}}(H, K^{-1}\Bbbk[x])$, define $f \in \Hom_{\cM^{H}}(H, \Bbbk[x])$ by
	\[ f(K^mx^n) = \begin{cases}
		\binom{n}{n-m}_{q^{2}} \frac{1}{(m+1)_{q^{2}}} \vep(g(K^mx^{m+1})) x^{n-m}, & 0 \leqslant m \leqslant n \\
		0, & \text{otherwise}.
	\end{cases} \]
	A direct calculation shows that $f$ is $H$-colinear and that $\vep(g(K^{-1})) \pi_{-1} + d^1 \circ f = g$.
	Thus, the class of $\pi_{-1}$ generates the cohomology group $\Ext_{\cM^{H}}^1(H, \Bbbk)$; that is, $\Ext_{\cM^{H}}^1(H, \Bbbk) = \Bbbk\,\overline{\pi_{-1}}$.
	Moreover, the injective dimension of the trivial $H$-comodule $\Bbbk$ is at least $1$.
	Therefore, $H$ is AS-co-regular of dimension $1$. Let $\int^{H}_r = \Ext_{\cM^{H}}^1(H, \Bbbk)$.
	
	For any $h^* \in H^*$, we have
	\[ (h^* \rhd \pi_{-1})(K^{-1}) = \langle h^*, K^{-1} \rangle \pi_{-1}(K^{-1}), \]
	so $h^* \rhd \pi_{-1} - \langle h^*, K^{-1} \rangle \pi_{-1} \in \mathrm{Im} (d^1_*)$.
	Thus, as a right $H$-comodule, $\textstyle \int^{H}_r = \Bbbk \overline{\pi_{-1}} \cong \Bbbk K^{-1}$.
	Let $\widetilde{\alpha}: H \to \Bbbk$ be a character such that $\textstyle \int^{H}_r \cong {}_{\widetilde{\alpha}}\Bbbk$.
	Since
	\[ (K \rhd \pi_{-1})(K^{-1}) = \pi_{-1}(K^{-1} K^{-1} K) \lhd K^{-1} = K^{-1} \lhd K^{-1} = q^2K^{-1}, \]
	it follows that $K \rhd \pi_{-1} - q^2 \pi_{-1} \in \mathrm{Im}(d^1_*)$, and hence $\widetilde{\alpha}(K) = q^2$.
	Moreover, $\widetilde{\alpha}(E) = \widetilde{\alpha}(KEK^{-1}) = \widetilde{\alpha}(q^2E)$.
	Since $q^2 \neq 1$, it follows that $\widetilde{\alpha}(E) = 0$.
  Thus $(\lint)^* \cong \rcoint$ as objects of $\YDlr(\id_H, S^4)$, where the twisted left-right Yetter--Drinfeld module structure on $(\lint)^*$ is given by \Cref{lem:dual-YD}.
\end{ex}

\subsection{Bijectivity of the antipode}\label{subsec:bijective}

As a first consequence of the $S^4$-formulas, we consider the bijectivity of the antipode.
Let $H$ be a Hopf algebra over a field $\Bbbk$ with antipode $S$.
The antipode is known to be bijective in several settings.

\begin{rem}\label{rem:bijective-known}
Let $H$ be a Hopf algebra with antipode $S$. Then $S$ is bijective in each of the following cases:
\begin{enumerate}[label=\textup{(\arabic*)},leftmargin=2.2em]
	\item \cite[Proposition~2]{LS69} $H$ is finite-dimensional;
	\item \cite[Proposition~2]{Rad77c} $H$ is co-Frobenius, that is, $H^*$ contains a nonzero integral;
	\item \cite{Rad92} $(H, R)$ is a quasitriangular Hopf algebra; 
	\item \cite[Theorem~A]{Skr06} $H$ can be embedded into a left perfect ring $Q$ such that $Q$ is an essential extension of $H$ as a right $H$-module;
	\item \cite[Theorem~0.2]{LOWY18} $H$ is AS-Gorenstein and every right invertible element is regular;
	\item \cite[Theorem~6.1]{Sci26} $H$ is either right or left noetherian.
\end{enumerate}
In case \textup{(1)}, bijectivity also follows from \eqref{eq:Radford}. In case \textup{(2)}, it follows from a formula of the same shape. In both formulas, the map defined by the right-hand side is the composite of conjugation by a group-like element and the bijection $h\mapsto\alpha\rhu h\lhu\alpha^{-1}$, so this map is bijective. Hence $S^4$ is bijective, and therefore $S$ is bijective.
\end{rem}

For any Hopf algebra $H$, the bialgebra $H^{\op\cop}$, with the opposite multiplication and the opposite comultiplication, is a Hopf algebra with antipode $S$.
Consequently, $H$ is left AS (co-)Gorenstein if and only if $H^{\op\cop}$ is right AS (co-)Gorenstein.
The following proposition follows from the $S^4$-formulas established in \Cref{thm:main} and \Cref{thm:main-co}.

\begin{prop}\label{prop:bijective-antipode-ASG}
	\begin{enumerate}[label=\textup{(\arabic*)},leftmargin=2.2em]
		\item Let $H$ be a left AS-Gorenstein Hopf algebra of dimension $d$.
		If the coradical of $H$ is coreflexive, or if ${}_{\vep}\Bbbk$ is of type $\mathrm{FP}_{d}$, then the antipode of $H$ is bijective.
		\item Let $H$ be a right AS-co-Gorenstein Hopf algebra.
		If the coradical of $H$ is coreflexive, then the antipode of $H$ is bijective.
	\end{enumerate}
\end{prop}
\begin{proof}
(1)
By \Cref{thm:main},
\[
S^4(h)=g\,(\alpha\rhu h\lhu\alpha^{-1})\,g^{-1}
\qquad\text{for all }h\in H.
\]
The maps $h\mapsto ghg^{-1}$ and $h\mapsto\alpha\rhu h\lhu\alpha^{-1}$ are bijections, so $S^4$, and hence $S$, is bijective.

(2) The same argument, applied to \eqref{eq:main-co-formula}, shows that $S$ is bijective.
\end{proof}

Any left AS-Gorenstein Hopf algebra $H$ satisfying the condition in \Cref{prop:bijective-antipode-ASG}(1) has a bijective antipode $S$. Since $S$ induces an algebra isomorphism $H \cong H^{\op}$ that interchanges the left and right AS-Gorenstein conditions, $H$ is also right AS-Gorenstein and hence AS-Gorenstein.
Any right AS-co-Gorenstein Hopf algebra $H$ satisfying the condition in \Cref{prop:bijective-antipode-ASG}(2) has a bijective antipode $S$. Since $S$ induces a coalgebra isomorphism $H \cong H^{\cop}$ that interchanges the left and right AS-co-Gorenstein conditions, $H$ is also left AS-co-Gorenstein and hence AS-co-Gorenstein.

\begin{rem}
	Let $H$ be an AS-Gorenstein Hopf algebra of dimension $d$.
	Suppose that the coradical of $H$ is coreflexive, or that ${}_{\vep}\Bbbk$ is of type $\mathrm{FP}_{d}$.
    \begin{enumerate}[label=\textup{(\arabic*)},leftmargin=2.2em]
    \item Let $K=H^{\op\cop}$. The antipode gives a Hopf isomorphism $S\colon H\to K$. As algebras, $K$ coincides with $H^{\op}$, so $\rint=\Ext^d_{H^{\op}}(\Bbbk,H)$ is the space of left homological integrals of $K$. By \Cref{thm:main}, $\lint$ is an $(\id_H,S^4)$-twisted right-left Yetter--Drinfeld module. Write $t\lhd h$ for its right $H$-action and $t\mapsto t_{-1}\otimes t_0$ for its left $H$-coaction. Transport along $S$ gives the right $K$-action and left $K$-coaction
    \[
    t\lhd_K S(h)=t\lhd h, \qquad t\mapsto S(t_{-1})\otimes t_0
    \]
    on $\rint$. A right $K$-action is a left $H$-action and a left $K$-coaction is a right $H$-coaction, so these become
    \[
    h\rhd t=t\lhd S^{-1}(h), \qquad t\mapsto t_0\otimes S(t_{-1}).
    \]
    Under this change of sides, $\rint$ is an $(S^4,\id_H)$-twisted left-right Yetter--Drinfeld module over $H$. If $\lint\cong\Bbbk_\alpha$ as a right $H$-module, then $\rint\cong{}_{\alpha^{-1}}\Bbbk$ as a left $H$-module. If $\lint\cong\Bbbk g$ as a left $H$-comodule, then $\rint\cong\Bbbk g^{-1}$ as a right $H$-comodule.
	
    \item The space of left homological integrals $\lint$ can be viewed as a right $H^{\circ} \#_{\id_H, S^4} H$-module.
	The formula \eqref{eq:main-formula} for $S^4$ can then be written as
	\[ S^4(h) = (\alpha \# g) (1 \# h) (\alpha \# g)^{-1} = (\alpha \# g) (1 \# h) (\alpha^{-1} \# g^{-1}), \]
	where $H$ is identified with a subalgebra of $H^{\circ} \#_{\id_H, S^4} H$ via $h \mapsto 1 \# h$.
	\end{enumerate}
\end{rem}


\subsection{\texorpdfstring{Order of $S^2$}{Order of S\texttwosuperior}}\label{subsec:order}

In this subsection, we apply the $S^4$-formulas to investigate when the antipode has finite order. The following result is due to Brown and Zhang.

\begin{prop}\cite[Proposition~6.2]{BZ08}\label{prop:inner}
	Let $H$ be a noetherian affine PI Hopf algebra.
	If $H$ has finite global dimension or is a finite module over its center, then some power of the antipode $S$ of $H$ is an inner automorphism.
\end{prop}
 
In the same paper, Brown and Zhang posed the following question. 

\begin{ques}\cite[Question~6.2]{BZ08}\label{q:finite-order}
	If $H$ is a noetherian affine PI Hopf algebra, does the antipode of $H$ have finite order?
\end{ques}

The $S^4$-formula provides a partial answer to \Cref{q:finite-order}.
We first outline the proof strategy.
By the $S^4$-formula (\Cref{thm:main} or \Cref{thm:main-co}), there exist a group-like element $g \in G(H)$ and a character $\alpha \colon H \to \Bbbk$ such that
\[ S^4(h) = g\, (\alpha \rhu h \lhu \alpha^{-1})\, g^{-1} \quad \text{for all } h \in H. \]
Let $\ad_g, \ad_{\alpha} \in \Aut_{\mathrm{Hopf}}(H)$ denote the automorphisms
\[ \ad_g(h) = g h g^{-1}, \quad \ad_{\alpha}(h) = \alpha \rhu h \lhu \alpha^{-1}, \]
respectively.
Then the $S^4$-formula can be rewritten as
\[ S^4 = \ad_g \circ \ad_{\alpha}. \]
Since $\ad_g$ and $\ad_{\alpha}$ commute, $S^4$ has finite order whenever both $\ad_g$ and $\ad_{\alpha}$ have finite order.
Thus, to prove that $S^2$ has finite order, it suffices to show that both $\ad_g$ and $\ad_{\alpha}$ have finite order.

Let $B$ be a right coideal subalgebra of $H$, and write $B^+ = B \cap \Ker \vep$.
If $B$ is central in $H$, then $HB^+ = B^+H$, so $HB^+$ is a two-sided ideal of $H$.
Moreover, $HB^+$ is a coideal: $\vep(HB^+)=0$, and $\Delta(HB^+)\subseteq HB^+\otimes H+H\otimes HB^+$.
Since $\Delta(B) \subseteq B \otimes H$, $S(b) = (\vep(b_{1}) - b_1)S(b_2) \in B^+H$ for $b\in B^{+}$.
Consequently, $S(HB^+) \subseteq B^+H = HB^+$.
Thus, $HB^+$ is stable under the antipode and is therefore a Hopf ideal of $H$.
Finally, if $H$ is finitely generated as a $B$-module, then $H/HB^+$ is finitely generated over $B/B^+ = \Bbbk$ and is therefore finite-dimensional.

\begin{prop}\label{prop:order-ASG-B}
	Let $H$ be an AS-Gorenstein Hopf algebra of dimension $d$ with antipode $S$.
	Suppose that the following conditions hold.
	\begin{enumerate}[label=\textup{(\arabic*)},leftmargin=2.2em]
		\item $H$ is finite as a module over a central right coideal subalgebra $B$.
		\item $H$ is faithfully flat as a left $B$-module.
		\item The coradical of $H$ is coreflexive, or ${}_{\vep}\Bbbk$ is of type $\mathrm{FP}_{d}$.
	\end{enumerate}
	Then $S^2$ has finite order, which divides $2 \dim_{\Bbbk} \overline{H}$, where $\overline{H} := H/HB^+$.
\end{prop}
\begin{proof}
	Let $\pi \colon H \to \overline{H}$, $h \mapsto \overline{h}$, be the canonical quotient Hopf algebra homomorphism.
	Set $n := \dim_{\Bbbk} \overline{H}$.
	Let $g \in G(H)$ be the group-like element and $\alpha\colon H\to\Bbbk$ the character afforded by \Cref{thm:main}.
	Then \eqref{eq:main-formula} can be rewritten as
	\[ S^4 = \ad_g \circ \ad_{\alpha} = \ad_{\alpha} \circ \ad_g. \]
	
	Since a finite-dimensional Hopf algebra is  free over any of its Hopf subalgebras \cite[Theorem~3.1.5]{Mon93}, the order of $G(\overline{H})$ divides $n$.
	The image $\overline{g}$ of $g$ in $\overline{H}$ is a group-like element, so $\overline{g}^{n} = \overline{1_H}$; that is, $g^{n} \in 1_H + HB^+$.
	Then
	\[ (\pi \otimes \id_H) \Delta(g^{n}) = \overline{g^{n}} \otimes g^{n} = \overline{1_H} \otimes g^{n}, \]
	so $g^{n} \in{}^{\mathrm{co}\overline{H}}H$.
	Since $H$ is faithfully flat over $B$, Takeuchi's correspondence \cite[Theorem~1]{Tak79} implies that ${}^{\mathrm{co}\overline{H}}H = B$.
	Hence $g^{n} \in B$, which is contained in the center $Z(H)$ of $H$.
	Therefore $\ad_{g^{n}} = \id_H$.
	
	For $z \in Z(H)$, right multiplication by $z$ on $H$ coincides with left multiplication by $z$, which is $H$-linear and, by the argument of \Cref{prop:trivial-ASG}(1), induces $\vep(z) \id$ on $\lint$.
	In particular, $\alpha(B^+) = \vep(B^+) = 0$, so $\alpha$ factors through an algebra homomorphism $\overline{\alpha} \colon \overline{H} \to \Bbbk$, which is a group-like element of the finite-dimensional Hopf algebra $\overline{H}^*$.
	Since the order of $G(\overline{H}^*)$ divides $n$, $\overline{\alpha}^{n} = \vep$; equivalently, $\alpha^{n} = \vep$ on $H$.
	
	Since $\ad_{g^{n}} = \id_H$ and $\alpha^n = \vep$, it follows that
	\[ S^{4n} = (\ad_g \circ \ad_{\alpha})^n = \ad_g^n \circ \ad_{\alpha}^n = \ad_{g^n} \circ \ad_{\alpha^n} = \id. \]
	Thus, $S^2$ has finite order dividing $2n$.
\end{proof}

\begin{rem}
	In the proof of \Cref{prop:order-ASG-B}, the order of $\alpha \in G(H^{\circ})$ is finite under the weaker assumption that $H$ is an affine Hopf algebra, finite as a module over its center \cite[Proposition~6.2(a) and the sentence following its proof]{BZ08}.
\end{rem}

To apply \Cref{prop:order-ASG-B} to Hopf algebras that are finite over their central Hopf subalgebras, we must verify the faithful flatness assumption.
The following lemma establishes that this assumption holds.

\begin{lem}\label{lem:weakly-finite}
	Let $H$ be a Hopf algebra that is finitely generated as a module over its center. Then $H$ is faithfully flat over any of its commutative Hopf subalgebras.
\end{lem}
\begin{proof}
    By \cite[Theorem 5.4]{Skr25}, a weakly finite Hopf algebra is faithfully flat over any of its commutative Hopf subalgebras. Thus it remains only to prove that $H$ is weakly finite.
	Recall that a ring $R$ is Dedekind-finite if $xy = 1$ implies $yx = 1$ for all $x, y \in R$, and that $R$ is weakly finite if the matrix ring $M_n(R)$ is Dedekind-finite for every $n \geqslant 1$.
	Since $H$ is finitely generated as a $Z(H)$-module, the canonical map $H \hookrightarrow \prod\limits_{\mathfrak{m} \in \Max Z(H)} H_{\mathfrak{m}}$ is injective.
	Moreover, each localization $H_{\mathfrak{m}}$ is a finite module over the central local ring $Z(H)_{\mathfrak{m}}$ and is therefore semilocal.
	Matrix rings over a semilocal ring are semilocal (see \cite[20.4]{Lam01}), and every semilocal ring is Dedekind-finite by \cite[Proposition~20.8]{Lam01}; thus, semilocal rings are weakly finite.
	Since weak finiteness is preserved under subrings and arbitrary products (see \cite[Corollary~1.10 and Proposition~1.9]{Lam99}), it follows that $H$ is weakly finite.
\end{proof}

Applying \Cref{prop:order-ASG-B} yields the following corollary, with faithful flatness guaranteed by \Cref{lem:weakly-finite}.

\begin{cor}\label{cor:order-ASG}
	Let $H$ be an AS-Gorenstein Hopf algebra of dimension $d$.
	Suppose that either the coradical of $H$ is coreflexive or ${}_{\vep}\Bbbk$ is of type $\mathrm{FP}_{d}$.
	If $H$ is a finite module over a central Hopf subalgebra $B$, then $S^2$ has finite order dividing $2\dim \overline{H}$, where $\overline{H} := H/HB^+$.
\end{cor}

\begin{prop}\label{prop:order-AScoG-B}
	Let $H$ be an AS-co-Gorenstein Hopf algebra with antipode $S$, and 
	suppose that the coradical of $H$ is coreflexive.
	If $H$ is a finite module over a central Hopf subalgebra $B$, then $S^2$ has finite order dividing $2\dim \overline{H}$, where $\overline{H} := H/HB^+$.
\end{prop}
\begin{proof}
	Let $\pi \colon H \to \overline{H}$, $h \mapsto \overline{h}$, denote the canonical quotient Hopf algebra homomorphism.
	Set $n := \dim_{\Bbbk}\overline{H}$.
	Let $\widetilde{g} \in G(H)$ be the group-like element and $\widetilde{\alpha}\colon H\to\Bbbk$ the character afforded by \Cref{thm:main-co}.
	Then \eqref{eq:main-co-formula} can be rewritten as
	\[ S^4 = \ad_{\widetilde{g}} \circ \ad_{\widetilde{\alpha}} = \ad_{\widetilde{\alpha}} \circ \ad_{\widetilde{g}}. \]
	
	By \Cref{lem:weakly-finite}, $H$ is faithfully flat as a left $B$-module.
	The same argument as in the proof of \Cref{prop:order-ASG-B} then shows that $\widetilde{g}^{n} \in B$, and hence $\ad_{\widetilde{g}^{n}} = \id_H$.
	By \Cref{prop:coresolution}, there exists an exact sequence
	\[ 0 \To \Bbbk \stackrel{d^0}{\To} M^0 \boxplus H \stackrel{d^1}{\To} M^1 \boxplus H \stackrel{d^2}{\To} \cdots \stackrel{d^{i}}{\To} M^i \boxplus H \stackrel{d^{i+1}}{\To} \cdots \]
	in $\YDrr$.
	For any right $H$-module $M$ on which $B$ acts via the counit, that is, $m \lhd b = \vep(b)m$, the right $B$-action on $M \boxplus H$ is given by
	\[ (m \boxplus h) \lhd b = m \lhd b_2 \boxplus S(b_1) h b_3 = m \boxplus S(b_1) h b_2 = m \boxplus h S(b_1) b_2 = \vep(b) (m \boxplus h), \]
	where the third equality follows from the centrality of $B$ in $H$.
	Hence, the $B$-action on each $M^i \boxplus H$ is induced by the counit.
	Moreover, for any $f \in \Hom_{\cM^H}(H, M^i \boxplus H)$, the action \eqref{eq:co-Hom-H-action} gives
	\[ (b \rhd f)(h) = f\bigl( S^3(b_3) h b_1 \bigr) \lhd S(b_2) = f\bigl( S^3(b_2) h b_1 \bigr) =  f\bigl( h b_1 S^3(b_2) \bigr) = f\bigl( h b_1 S(b_2) \bigr) = \vep(b) f(h), \]
	which shows that the $B$-action on each $\Hom_{\cM^H}(H, M^i \boxplus H)$ is also induced by the counit $\vep$.
	Therefore, the $B$-action on $\rcoint$ factors through the counit. 
	Hence $\widetilde{\alpha}|_B = \vep|_B$, so $\widetilde{\alpha}(HB^+) = 0$, and $\widetilde{\alpha}$ factors through $\pi$ to give an element of $G(\overline{H}^*)$.
	Since $\overline{H}$ is finite-dimensional, the characters of $\overline{H}$ have finite order dividing $\dim \overline{H} = n$; thus $\widetilde{\alpha}^{n} = \vep$.
	
	Since $\ad_{\widetilde{g}^{n}} = \id_H$ and $\widetilde{\alpha}^n = \vep$, it follows that
	\[ S^{4n} = (\ad_{\widetilde{g}} \circ \ad_{\widetilde{\alpha}})^n = \ad_{\widetilde{g}}^n \circ \ad_{\widetilde{\alpha}}^n = \ad_{\widetilde{g}^n} \circ \ad_{\widetilde{\alpha}^n} = \id. \]
	Hence, $S^2$ has finite order dividing $2n$.
\end{proof}

Two analogues of \Cref{prop:order-ASG-B} hold when finite generation and faithful flatness over a central right coideal subalgebra are replaced by finite cogeneration and faithful coflatness over a cocentral quotient Hopf algebra: \Cref{prop:order-AScoG-C} in the AS-co-Gorenstein case, and \Cref{prop:order-ASG-C} in the AS-Gorenstein case.
We now recall the relevant definitions.
Let $C$ be a coalgebra, and let $N$ be a left $C$-comodule.
The comodule $N$ is a \emph{finitely cogenerated} $C$-comodule if it embeds, as a $C$-comodule, into a finite direct sum of copies of $C$.
A coalgebra map $\pi \colon C \to D$ is called \emph{cocentral} if, for all $c \in C$,
\[ \pi(c_1) \otimes c_2 = \pi(c_2) \otimes c_1. \]

Let $\pi \colon H \to C$ be a quotient map onto a left $H$-module coalgebra $C$, and let
\[ B := {}^{\mathrm{co}\pi}H = \{ h \in H \mid (\pi \otimes \id) \Delta(h) = \pi(1_H) \otimes h\} \]
denote the subalgebra of left $C$-coinvariants.
Then $B$ is a right coideal of $H$.
If $\pi$ is cocentral, then ${}^{\mathrm{co}\pi}H = H^{\mathrm{co}\pi} := \{ h \in H \mid (\id \otimes \pi) \Delta(h) = h \otimes \pi(1_H) \}$.
Thus, $B$ is also a left coideal and a sub-bialgebra of $H$.
Moreover, for $b \in B$, we have
\begin{align*}
	b \otimes \pi(1_H) & = b_3 \otimes \pi\bigl( S(b_1) b_2 \bigr) \\
	& = b_3 \otimes S(b_1) \pi(b_2) & \text{ ($\pi$ is $H$-linear)} \\
	& = b_2 \otimes S(b_1) \pi(1_H) & \text{ (since $B$ is a sub-bialgebra and $B = {}^{\mathrm{co}\pi}H$)} \\
	& = b_2 \otimes \pi\bigl( S(b_1) \bigr) & \text{ ($\pi$ is $H$-linear)}
\end{align*}
Applying $S \otimes \id_C$ to the above equality yields $S(b) \otimes \pi(1_H) = S(b_2) \otimes \pi(S(b_1))$; equivalently, $(\id \otimes \pi)\Delta(S(b)) = S(b) \otimes \pi(1_H)$.
Hence $S(B) \subseteq H^{\mathrm{co}\pi} = B$, so $B$ is a Hopf subalgebra of $H$.
If $H$ embeds into $C^{\oplus r}$ as a left $C$-comodule for some $r \geqslant 1$, then $B$ embeds into ${}^{\mathrm{co}C}(C^{\oplus r}) = ({}^{\mathrm{co}C}C)^{\oplus r} = \Bbbk^r$ and is therefore finite-dimensional.

\begin{prop}\label{prop:order-AScoG-C}
	Let $H$ be an AS-co-Gorenstein Hopf algebra of dimension $d$ with antipode $S$, and let $\pi \colon H \to C$ be a quotient map onto a left $H$-module coalgebra $C$.
	Suppose that the following conditions hold.
	\begin{enumerate}[label=\textup{(\arabic*)},leftmargin=2.2em]
		\item $H$ is finitely cogenerated and faithfully coflat as a left $C$-comodule.
		\item $\pi$ is cocentral.
		\item The coradical of $H$ is coreflexive.
	\end{enumerate}
	Then $S^2$ has finite order dividing $2 \dim_{\Bbbk} B$, where $B = {}^{\mathrm{co}\pi}H$.
\end{prop}
\begin{proof}
	Set $n :=\dim_{\Bbbk}B$.
	Let $\widetilde{g} \in G(H)$ be the group-like element and $\widetilde{\alpha} \colon H \to \Bbbk$ the character afforded by \Cref{thm:main-co}.
	Then \eqref{eq:main-co-formula} can be rewritten as
	\[ S^4 = \ad_{\widetilde{g}} \circ \ad_{\widetilde{\alpha}} = \ad_{\widetilde{\alpha}} \circ \ad_{\widetilde{g}}. \]
	
	Applying \Cref{prop:coresolution} to $V = \Bbbk_{\vep}^{\eta} \in \YDrr$ yields an exact sequence
	\[ 0 \To \Bbbk_{\vep}^{\eta} \stackrel{d^0}{\To} G^0 \stackrel{d^1}{\To} G^1 \stackrel{d^2}{\To} \cdots \stackrel{d^{i}}{\To} G^i \stackrel{d^{i+1}}{\To} \cdots \]
	in $\YDrr$, in which each $G^i$ is an injective $H$-comodule.
	We identify $C^*$ with its image in $H^*$ under $\pi^*$.
	For any $\lambda \in \Hom_{\cM^H}(H, G^i)$ and $c^* \in C^*$, the action \eqref{eq:co-Hom-H^*-action} yields
	\[ \bigl( (c^*\pi) \rhd \lambda \bigr)(h) = \langle c^*, \pi(h_1) \rangle \lambda(h_2) = \langle c^*, \pi(h_2) \rangle \lambda(h_1) = \langle c^*, \pi(\lambda(h)_1) \rangle \lambda(h)_0 = (c^* \pi) \rhd \lambda(h). \]
	Thus, the $C^*$-action \eqref{eq:co-Hom-H^*-action} on each $\Hom_{\cM^H}(H, G^i)$ coincides with the action induced by the $C^*$-action on $G^i$.
	Since $\pi$ is cocentral, the action map $c^*\pi \rhd -$ is right $H$-colinear on each $G^i$ and thus defines an $H$-colinear chain endomorphism $\xi$ of $G^{\bullet}$.
	Moreover, since $c^*\pi \rhd -$ coincides with $\langle c^*, \pi(1_H) \rangle \id$ on $\Bbbk_{\vep}^{\eta}$, $\xi \colon G^{\bullet} \to G^{\bullet}$ is homotopic to $\langle c^*, \pi(1_H) \rangle \id_{G^{\bullet}}$.
	Applying $\Hom_{\cM^H}(H, -)$, we conclude that the induced chain endomorphism $\Hom_{\cM^H}(H, \xi)$ of $\Hom_{\cM^H}(H, G^{\bullet})$ is homotopic to $\langle c^*, \pi(1_H) \rangle \id$.
	Since the chain endomorphism induced by the action of $c^* \pi$ on $\Hom_{\cM^H}(H, G^{\bullet})$ coincides with $\Hom_{\cM^H}(H, \xi)$, it is also homotopic to $\langle c^*, \pi(1_H) \rangle \id$.
	Hence, the induced action of $c^*\pi$ on $\Ext_{\cM^H}^i(H, \Bbbk^{\eta})$ is given by the scalar $\langle c^*, \pi(1_H) \rangle$.
	In particular, fix $0 \neq x \in \rcoint$ and $c^* \in C^*$.
	Then
	\[ \langle c^*, \pi(1_H) \rangle x = (c^*\pi) \rhd x = \langle c^*, \pi(x_{1}) \rangle x_0 = \langle c^*, \pi(\widetilde{g}) \rangle x. \]
	Since $x \neq 0$ and $c^*$ is arbitrary, it follows that $\pi(1_H) = \pi(\widetilde{g})$, and hence $\widetilde{g} \in B$.
	Since $B$ is finite-dimensional, its group-like elements have finite order; in particular, the order of $\widetilde{g}$ divides $n$.
	Consequently, $\widetilde{g}^{n} = 1_H$.
	
	The restriction of $\widetilde{\alpha}$ to the finite-dimensional Hopf subalgebra $B$ has finite order dividing the dimension of $B$; hence $\widetilde{\alpha}^{n}|_{B} = \vep|_{B}$.
	It follows that $\widetilde{\alpha}^{n}(HB^+) = \widetilde{\alpha}^{n}(H)\widetilde{\alpha}^{n}(B^+) = \widetilde{\alpha}^{n}(H)\vep(B^+) = 0$.
	Since $H$ is a faithfully coflat $C$-comodule, Takeuchi's correspondence \cite[Theorem~2]{Tak79} yields $\Ker(\pi) = HB^+$.
	Hence $\widetilde{\alpha}^{n}$ is cocentral because it factors through the cocentral map $\pi$.
	Thus, for any $h \in H$, we have
	\[ \ad_{\widetilde{\alpha}^{n}}(h) = \langle \widetilde{\alpha}^{n}, S(h_1) \rangle h_2 \langle \widetilde{\alpha}^{n}, h_3 \rangle = \langle \widetilde{\alpha}^{n}, S(h_1) \rangle \langle \widetilde{\alpha}^{n}, h_2 \rangle h_3 = h. \]
	It follows that $\ad_{\widetilde{\alpha}^{n}} = \id_H$.
	
	Since $\widetilde{g}^n = 1_H$, $\ad_{\widetilde{\alpha}^{n}} = \id_H$, and $\ad_{{\widetilde{g}}}$ commutes with $\ad_{\widetilde{\alpha}}$, we obtain
	\[ S^{4n} = (\ad_{\widetilde{g}} \circ \ad_{\widetilde{\alpha}})^n = \ad_{\widetilde{g}}^n \circ \ad_{\widetilde{\alpha}}^n = \ad_{\widetilde{g}^n} \circ \ad_{\widetilde{\alpha}^n} = \id_H. \]
	Hence $S^2$ has finite order dividing $2n$.
\end{proof}

\begin{prop}\label{prop:order-ASG-C}
	Let $H$ be an AS-Gorenstein Hopf algebra of dimension $d$ with antipode $S$, and let $\pi \colon H \to K$ be a surjective Hopf algebra homomorphism.
	Suppose that the following conditions hold.
	\begin{enumerate}[label=\textup{(\arabic*)},leftmargin=2.2em]
		\item $H$ is finitely cogenerated and faithfully coflat as a left $K$-comodule.
		\item $\pi$ is cocentral.
		\item The coradical of $H$ is coreflexive, or ${}_{\vep}\Bbbk$ is of type $\mathrm{FP}_{d}$.
	\end{enumerate}
	Then $S^2$ has finite order dividing $2 \dim_{\Bbbk} B$, where $B = {}^{\mathrm{co}\pi}H$.
\end{prop}
\begin{proof}
	Set $n :=\dim_{\Bbbk}B$.
	Let $g \in G(H)$ be the group-like element and $\alpha \colon H \to \Bbbk$ be the character afforded by \Cref{thm:main}.
	Then \eqref{eq:main-formula} can be rewritten as
	\[ S^4 = \ad_{g} \circ \ad_{\alpha} = \ad_{\alpha} \circ \ad_{g}. \]
	
	For any left $H$-comodule $N$ whose induced $K$-coaction is given by the group-like element $\pi(1_H)$, the coaction on $H \boxtimes N$ satisfies
	\begin{align*}
		(\pi \otimes \id) \rho(h \boxtimes n) & = \pi\bigl( h_1 n_{-1} S(h_3) \bigr) \otimes h_2 \boxtimes n_0 \\
		& = \pi(h_1) \pi(n_{-1}) \pi\bigl( S(h_3) \bigr) \otimes h_2 \boxtimes n_0 \\
		& = \pi(h_1) \pi(1_H) S\bigl( \pi(h_3) \bigr) \otimes h_2 \boxtimes n \\
		& = \pi(h_1) \pi(1_H) S\bigl( \pi(h_2) \bigr) \otimes h_3 \boxtimes n & \text{ ($\pi$ is cocentral)} \\
		& = \pi\bigl( h_1 S(h_2) \bigr) \otimes h_3 \boxtimes n \\
		& = \pi(1_H) \otimes h \boxtimes n.
	\end{align*}
	Thus, the $K$-coaction on $H \boxtimes N$ is also induced by $\pi(1_H)$.
	Consider the exact sequence
	\[ \cdots \To H \boxtimes N_i \To \cdots \To H \boxtimes N_1 \To H \boxtimes N_0 \To \Bbbk \To 0 \]
	in \Cref{prop:resolution}, where the $K$-coaction on each $N_i$ is induced by $\pi(1_H)$.
	For any $h \in H$, $k^* \in K^*$, $n \in N_i$, and $f \in \Hom_H(H \boxtimes N_i, H)$, we have
	\begin{align*}
		&  \bigl( f \lhd (k^*\pi) \bigr)(h \boxtimes n) \\
		& \xlongequal{\eqref{eq:Hom-H^*-action}} f\bigl( (h \boxtimes n)_0 \bigr)_2 \bigl\langle k^* \pi, S^3\bigl( f\bigl( (h \boxtimes n)_0 \bigr)_{3} \bigr) S\bigl( (h \boxtimes n)_{-1} \bigr) f\bigl( (h \boxtimes n)_0 \bigr)_{1} \bigr\rangle \\
		& = f\bigl( (h \boxtimes n)_0 \bigr)_2 \bigl\langle k^*, \pi\Bigl( S^3\bigl( f\bigl( (h \boxtimes n)_0 \bigr)_{3} \bigr) \Bigr) \pi\Bigl( S\bigl( (h \boxtimes n)_{-1} \bigr) \Bigr) \pi \Bigl( f\bigl( (h \boxtimes n)_0 \bigr)_{1} \Bigr) \bigr\rangle \\
		& = f(h \boxtimes n)_2 \bigl\langle k^*, S^3 \Bigl( \pi \bigl( f(h \boxtimes n)_{3} \bigr) \Bigr) \pi \bigl( f(h \boxtimes n)_{1} \bigr) \bigr\rangle \\
		& = f(h \boxtimes n)_3 \bigl\langle k^*, S^3 \Bigl( \pi \bigl( f(h \boxtimes n)_{2} \bigr) \Bigr) \pi \bigl( f(h \boxtimes n)_{1} \bigr) \bigr\rangle \hspace{3.8cm} \text{ ($\pi$ is cocentral)} \\
		& = \langle k^*, \pi(1_H) \rangle f(h \boxtimes n),
	\end{align*}
	where the last equality holds because the surjectivity and cocentrality of $\pi$ imply that $K$ is cocommutative.
	Thus, the $K^*$-action on $\Hom_H(H \boxtimes N_i, H)$ factors through evaluation at $\pi(1_H)$.
	Therefore, the $K^*$-action on $\lint$ factors through evaluation at $\pi(1_H)$; that is, for any $x \in \lint$ and $k^* \in K^*$,
	\[ \langle k^*, \pi(1_H) \rangle x = x \lhd (k^*\pi) = \langle k^*, \pi(x_{-1}) \rangle x_0 = \langle k^*, \pi(g) \rangle x. \]
	Since $x \neq 0$ and $k^*$ is arbitrary, we obtain $\pi(1_H) = \pi(g)$, and hence $g \in B$.
	Since $B$ is finite-dimensional, its group-like elements have finite order; in particular, the order of $g$ divides $n$.
	Consequently, $g^{n} = 1_H$.
	
    Since $\alpha\colon H\to\Bbbk$ is a character, its restriction to the finite-dimensional Hopf subalgebra $B$ has finite order dividing $n$, so $\alpha^{n}|_{B}=\vep|_{B}$.
    Thus $\alpha^{n}(HB^{+})=0$.
    Faithful coflatness and Takeuchi's correspondence \cite[Theorem~2]{Tak79} give $\Ker(\pi)=HB^{+}$, as in the proof of \Cref{prop:order-AScoG-C}.
    Hence $\alpha^{n}$ factors through the cocentral map $\pi$, and is therefore cocentral.
    For any $h\in H$,
    \[
    \ad_{\alpha^{n}}(h)
    =\langle\alpha^{n},S(h_1)\rangle h_2\langle\alpha^{n},h_3\rangle
    =\langle\alpha^{n},S(h_1)\rangle\langle\alpha^{n},h_2\rangle h_3
    =h,
    \]
    so $\ad_{\alpha^{n}}=\id_H$.
	
	Since $g^n = 1_H$, $\ad_{\alpha^{n}} = \id_H$, and $\ad_g$ commutes with $\ad_{\alpha}$, we obtain
	\[ S^{4n} = (\ad_{g} \circ \ad_{\alpha})^n = \ad_{g}^n \circ \ad_{\alpha}^n = \ad_{g^n} \circ \ad_{\alpha^n} = \id_H. \]
	Hence $S^2$ has finite order dividing $2n$.
\end{proof}

\section*{Acknowledgments}

AI tools were used for language editing and grammar checking.
The authors take full responsibility for the manuscript.
The authors would like to thank Liyu Liu and Yinhuo Zhang for their valuable comments and suggestions.

\bibliographystyle{siam}
\bibliography{References}

\end{document}